\documentclass[a4paper, final, 12pt]{article}
\usepackage{amsmath, amssymb, latexsym, amscd, amsthm,amsfonts,amstext}
\usepackage[mathscr]{eucal}
\usepackage{graphicx}
\usepackage{subfig}
\usepackage{float}
\usepackage{color}
\usepackage{hyperref}
\usepackage[utf8]{inputenc}
\usepackage[english]{babel}

\numberwithin{equation}{section}
\input{tcilatex}
\begin{document}

\title{Convergent numerical methods for parabolic equations with reversed
time via a new Carleman estimate}
\author{Michael V. Klibanov\thanks{%
Department of Mathematics and Statistics, University of North Carolina
Charlotte, Charlotte, NC, 28223, mklibanv@uncc.edu.} and Anatoly G. Yagola%
\thanks{%
Department of Mathematics, Faculty of Physics, Moscow State University,
Moscow, 119991, Russia, e-mail: yagola@physics.msu.ru}}
\date{}
\maketitle

\begin{abstract}
The key tool of this paper is a new Carleman estimate for an arbitrary
parabolic operator of the second order for the case of reversed time data.
This estimate works on an arbitrary time interval. On the other hand, the
previously known Carleman estimate for the reversed time case works only on
a sufficiently small time interval. First, a stability estimate is proven.
Next, the quasi-reversibility numerical method is proposed for an arbitrary
time interval for the linear case. This is unlike a sufficiently small time
interval in the previous work. The convergence rate for the
quasi-reversibility method is established. Finally, the quasilinear
parabolic equation with reversed time is considered. A weighted globally
strictly convex Tikhonov-like functional is constructed.\ The weight is the
Carleman Weight Function which is involved in that Carleman estimate. The
global convergence of the gradient projection method to the exact solution
is proved for this functional.
\end{abstract}



\section{Introduction}

\label{sec:1}

In this paper, we construct convergent numerical methods for linear and
quasilinear parabolic equations with reversed time. The key tool of the
convergence analysis is a new Carleman estimate. While the previously known
Carleman estimate for these problems works only on a sufficiently small time
interval, the one of this paper works on any finite time interval.

All functions below are real valued ones. Below $x=\left(
x_{1},x_{2},...,x_{n}\right) $ denotes points in $\mathbb{R}^{n}$ and $%
\nabla f=\left( f_{x_{1}},f_{x_{2}},...,f_{x_{n}}\right) $ for any
appropriate function $f\left( x\right) .$ Let $\Omega \subset \mathbb{R}^{n}$
be a bounded domain with a piecewise smooth boundary $\partial \Omega .$ Let 
$T>0$ and $\tau \in \left( 0,T\right) $ be two numbers. Denote 
\begin{equation}
Q_{T}=\Omega \times \left( 0,T\right) ,S_{T}=\partial \Omega \times \left(
0,T\right) ,Q_{T\tau }=\Omega \times \left( \tau ,T\right) .  \label{1.1}
\end{equation}%
For $i,j=1,...,n$, let functions $a_{ij}\left( x,t\right) $ be such that%
\begin{equation}
a_{ij}\left( x,t\right) =a_{ij}\left( x,t\right) \in C^{1}\left( \overline{Q}%
_{T}\right) ,  \label{1.2}
\end{equation}%
\begin{equation}
\mu _{1}|\xi |^{2}\leq \dsum\limits_{i,j=1}^{n}a_{ij}\left( x,t\right) \xi
_{i}\xi _{j}\leq \mu _{2}|\xi |^{2},\forall \left( x,t\right) \in \overline{Q%
}_{T},\forall \xi \in \mathbb{R}^{n},  \label{1.3}
\end{equation}%
\begin{equation}
\mu _{1},\mu _{2}=const.>0,\mu _{1}\leq \mu _{2}.  \label{1.4}
\end{equation}%
Introduce a uniform elliptic operator $L$ of the second order in the domain $%
Q_{T},$ 
\begin{equation}
Lu=\dsum\limits_{i,j=1}^{n}a_{ij}\left( x,t\right) u_{x_{i}x_{j}},\text{ }%
\left( x,t\right) \in Q_{T}.  \label{1.5}
\end{equation}%
Let the function $F\left( y\right) \in C^{1}\left( \mathbb{R}^{2n+2}\right)
. $ For this function, we assume that there exists a constant $\overline{C}=%
\overline{C}\left( F\right) >0$ depending only on $F$ such that%
\begin{equation}
\left\vert F\left( y_{1}\right) -F\left( y_{2}\right) \right\vert \leq 
\overline{C}\left\vert y_{1}-y_{2}\right\vert ,\forall y_{1},y_{2}\in 
\mathbb{R}^{2n+2}.  \label{1.10}
\end{equation}%
Consider the following quasilinear parabolic equation:%
\begin{equation}
u_{t}=Lu+F\left( \nabla u,u,x,t\right) ,\left( x,t\right) \in Q_{T}.
\label{1.6}
\end{equation}%
We impose the zero Dirichlet boundary condition on the function $u,$%
\begin{equation}
u\mid _{S_{T}}=0.  \label{1.07}
\end{equation}%
Since we work with the time reversed case, we assume that the function $%
u\left( x,t\right) $ is known at the final time $T$, 
\begin{equation}
u\left( x,T\right) =g\left( x\right) ,x\in \Omega .  \label{1.8}
\end{equation}%
Thus, we have obtained the following problem:

\textbf{Problem with Time Reversed Data.} \emph{Suppose that conditions (\ref%
{1.1})-(\ref{1.5}). Find a function }$u\in C^{2,1}\left( \overline{Q}%
_{T}\right) $\emph{\ satisfying conditions (\ref{1.6})-(\ref{1.8}).}

One of possible applications is in the case when a solid is heated and the
initial temperature is unknown. However, one can measure the temperature of
this solid at a final time. It is required to restore the temperature
distribution inside of this solid at all preceding times. Another
application of this problem, which was found recently, is in financial
mathematics, more precisely in a the problem of forecast of prices of stock
options using the Black-Scholes equation and real market data \cite{KKG}.

The problem of our interest is well known to be unstable, i.e. this is an
ill-posed problem. H\"{o}lder stability estimates for this problem are
known, see, e.g. \cite{Is,KAP,LRS}. To the best knowledge of the author, the
strongest H\"{o}lder stability result, which is valid for an arbitrary large
time interval $t\in \left( 0,T\right) ,$ is obtained by Isakov, see Theorem
3.1.3 in \cite{Is}. In section 2 of chapter 4 of \cite{LRS} and later in 
\cite{KAP} a Carleman estimate was used to obtain the H\"{o}lder stability
estimate. However, that H\"{o}lder stability estimate is valid only on a
sufficiently small time interval $t\in \left( T-\varepsilon ,T\right) $ for
a sufficiently small $\varepsilon >0.$ The smallness of this interval is due
to the Carleman Weight Function (CWF), which has been used in the Carleman
estimates for that problem so far \cite{KAP,LRS}. This function is $\left(
k+T-t\right) ^{-2\lambda },$ where numbers $k,T>0$ are sufficiently small
and the parameter $\lambda >0$ is sufficiently large. The same CWF was used
in \cite{Fr} for the proof of the uniqueness theorem.

Using this estimate, the first author has constructed in \cite{KAP} the
quasi-reversibility method (QRM) for the above problem in the linear case
and has proven its convergence, again on a sufficiently small time interval.
Below as well as in \cite{KAP,KKG} the QRM is realized via the minimization
of a certain regularization functional. On the other hand, QRM is quite
often realized via a proper perturbation of the underlying PDE operator \cite%
{Is,LL}. A surprising idea of the recent publication of Kaltenbacher and
Rundell \cite{Kalt} it to use the non local operator of the fractional $t-$%
derivative as the perturbing operator for the QRM. In the work of Tuan, Khoa
and Au \cite{Tuan} another version of the QRM is constructed for the
quasilinear case. Its convergence was proven for an arbitrary $T$. However,
the perturbation operator of \cite{Tuan} is a very complicated one. Another
version of the QRM was proposed in \cite{Is}.\ This version works only in
the case when the set of eigenfunctions of the underlying elliptic operator
forms an orthonormal basis in $L_{2}\left( \Omega \right) .$

The QRM was originally developed by Lattes and Lions in 1969 \cite{LL}.
Their idea became quite popular since then with many publications treating a
variety of ill-posed problems for PDEs.\ In this regard we refer to, e.g. 
\cite{Bour1,Bour2,ClasKlib,Is,KS,KKG,KN,Tuan}. In particular, the first
author has shown in the survey paper \cite{KAP} that as soon as a proper
Carleman estimate for an ill-posed problem for a linear PDE is available,
then the QRM can be constructed for this problem and its convergence rate
can be established. In the case of the time reversed data for the linear
parabolic PDE the construction of the QRM in \cite{KAP} is valid only on a
sufficiently small time interval $\left( 0,T\right) .$

As to the ill-posed problems for quasilinear PDEs, it was shown in \cite{KQ}
that, again as soon as a proper Carleman estimate is available for the
linear principal part of the PDE operator, a weighted globally strictly
convex Tikhonov-like functional can be constructed, i.e. this problem can be
\textquotedblleft convexified". The key element of this functional is the
presence of the CWF in it, i.e. the function which is involved as the weight
in the Carleman estimate for that linear principal part of the PDE operator.
In the follow up paper \cite{Bak} existence and uniqueness of the minimizer
of that functional were established and global convergence to the exact
solution of the gradient projection method of the minimization of this
functional was proven. As to the\ quasilinear parabolic equations, they were
considered in \cite{Bak,KQ,KKLY} only for the case of lateral Cauchy data
with numerical results in \cite{Bak,KKLY}. However, the case of time
reversed data was not considered in \cite{Bak,KQ,KKLY}.

The idea of convexifying coefficient inverse problems was first proposed in
1995-1997 in \cite{KI,Klib97}, although without numerical studies, also see 
\cite{KT}. Recently the interest to the convexification approach was
renewed. First, this was done only analytically \cite{BK1,KK}. Next, a
number of papers was published, which combined the theory with numerical
studies, see, e.g. \cite{Bak,KKLY} for the case of quasilinear PDEs and \cite%
{KlibKol1,KlibKol2,KEIT,Ktime} for coefficient inverse problems.

We call a numerical method for an ill-posed problem \emph{globally convergent%
} if there is a theorem claiming that it converges to the exact solution of
this problem without an a priori knowledge of a sufficiently small
neighborhood of this solution. On the other hand, we call a numerical method
for that problem \emph{locally convergent} if its convergence is rigorously
guaranteed only if iterations start in a sufficiently small neighborhood of
the exact solution. Thus, the convexification is a \emph{globally convergent}
method (see Theorem 5.4 in section 5). On the other hand, a gradient-like
method being applied to a non-convex Tikhonov-like functional, might have
guaranteed convergence to the exact solution only if its starting point is
located in a sufficiently small neighborhood of that solution. The latter is
called the \emph{local convergence}.

New elements of this paper are:

\begin{enumerate}
\item A new Carleman estimate for a general parabolic operator of the second
order with time reversed data is proven. This estimate works on an arbitrary
time interval $t\in \left( 0,T\right) $, unlike a sufficiently small
interval of previous publications \cite{Fr,KAP,LRS}. Results listed in items
2-4 below are based on this estimate.

\item A stability estimate is proven for the above Problem with Time
Reversed Data. This estimate is somewhat \textquotedblleft between" H\"{o}%
lder and logarithmic stability estimates. In other words, although it is
weaker than the H\"{o}lder stability estimate of \cite{Is,KAP,LRS}, it is
stronger than the logarithmic stability estimate. Still, the main advantage
of our stability estimate over ones in \cite{KAP,LRS} is that it works
without a smallness assumption imposed on the time interval.

\item In the linear case, the QRM is constructed, existence and uniqueness
of the minimizer as well as convergence of minimizers to the exact solution
are proven. Unlike previous works \cite{KAP,KKG}, a smallness assumption is
not imposed on the time interval.

\item In the quasilinear case, this problem is convexified \emph{for the
first time}. In other words, a weighted globally strictly convex
Tikhonov-like functional is constructed with the CWF in it. Both the
existence and uniqueness of its minimizer are proved.\ In addition, the
global convergence of the gradient projection method to the exact solution
is established. This way we avoid the use of a complicated perturbation
operator of \cite{Tuan}.
\end{enumerate}

In section 2 we prove the new Carleman estimate. In section 3 we present
stability estimate. In section 4 we describe the quasi-reversibility method
for the linear case, prove existence and uniqueness of the minimizer as well
as convergence of minimizers to the exact solution when the level of the
noise in the data tends to zero. In section 5 we construct the above
mentioned weighted globally strictly convex Tikhonov-like functional and
formulate corresponding theorems. These theorems are proved in section 6.

\section{Carleman Estimate}

\label{sec:2}

\textbf{Theorem 2.1}. (Carleman estimate). \emph{Assume that conditions (\ref%
{1.1})-(\ref{1.5}) are in place. Then there exists a sufficiently large
number }$\nu _{0}=\emph{\ }\nu _{0}\left( \mu _{1},\mu
_{2},\max_{i,j}\left\Vert a_{ij}\right\Vert _{C^{1}\left( \overline{Q}%
_{T}\right) },Q_{T}\right) >1$\emph{\ depending only on listed parameters
and the number }$C>0$\emph{\ depending on the same parameters as ones of }$%
\nu _{0}$\emph{, such that for all functions }$u\in H^{2}\left( Q_{T}\right) 
$\emph{\ satisfying the zero Dirichlet boundary condition (\ref{1.07}) the
following Carleman estimate is valid}%
\begin{equation*}
\dint\limits_{Q_{T}}\left( u_{t}-Lu\right) ^{2}\exp \left( 2\left(
t+1\right) ^{\nu }\right) dxdt
\end{equation*}%
\begin{equation}
\geq C\sqrt{\nu }\dint\limits_{Q_{T}}\left( \nabla u\right) ^{2}\exp \left(
2\left( t+1\right) ^{\nu }\right) dxdt+C\nu
^{2}\dint\limits_{Q_{T}}u^{2}\exp \left( 2\left( t+1\right) ^{\nu }\right)
dxdt  \label{1.9}
\end{equation}%
\begin{equation*}
-C\exp \left( 3\left( T+1\right) ^{\nu }\right) \left\Vert u\left(
x,T\right) \right\Vert _{L_{2}\left( \Omega \right) }^{2}-C\left\Vert \nabla
u\left( x,0\right) \right\Vert _{L_{2}\left( \Omega \right) }^{2},\text{ }%
\forall \nu \geq \nu _{0}.
\end{equation*}

\textbf{Proof}. In this proof, 
\begin{equation}
u\in C^{2}\left( \overline{Q}_{T}\right) ,u\mid _{S_{T}}=0.  \label{1}
\end{equation}%
The case $u\in H^{2}\left( Q_{T}\right) ,$ $u\mid _{S_{T}}=0$ can be
obtained from (\ref{1}) via density arguments. Everywhere below in this
paper $C=\emph{\ }C\left( \mu _{1},\mu _{2},\max_{i,j}\left\Vert
a_{ij}\right\Vert _{C^{1}\left( \overline{Q}_{T}\right) },Q_{T}\right) >0$
denotes different constants depending only on listed parameters.

We prove this theorem in six steps. Introduce a new function $v\left(
x,t\right) ,$ 
\begin{equation*}
v\left( x,t\right) =u\left( x,t\right) \exp \left( \left( t+1\right) ^{\nu
}\right) .
\end{equation*}%
Hence,%
\begin{equation*}
u=v\exp \left( -\left( t+1\right) ^{\nu }\right) ,
\end{equation*}%
\begin{equation*}
u_{t}=\left[ v_{t}-\nu \left( t+1\right) ^{\nu -1}v\right] \exp \left(
-\left( t+1\right) ^{\nu }\right) ,
\end{equation*}%
\begin{equation*}
u_{x_{i}x_{j}}=v_{x_{i}x_{j}}\exp \left( -\left( t+1\right) ^{\nu }\right) .
\end{equation*}%
Hence,%
\begin{equation*}
\left( u_{t}-Lu\right) ^{2}\exp \left( 2\left( t+1\right) ^{\nu }\right) = 
\left[ v_{t}-\left( Lv+\nu \left( t+1\right) ^{\nu -1}v\right) \right] ^{2}
\end{equation*}%
\begin{equation}
\geq v_{t}^{2}-2v_{t}\left( Lv+\nu \left( t+1\right) ^{\nu -1}v\right)
\label{2.1}
\end{equation}%
\begin{equation*}
=v_{t}^{2}-2v_{t}Lv-2\nu \left( t+1\right) ^{\nu -1}v_{t}v.
\end{equation*}

\textbf{Step 1}. First, we estimate from the below $-2v_{t}Lv,$%
\begin{equation}
-2v_{t}Lv=-\dsum\limits_{i,j=1}^{n}\left(
a_{i,j}v_{x_{i}x_{j}}v_{t}+a_{i,j}v_{x_{j}x_{i}}v_{t}\right) .  \label{2.2}
\end{equation}%
Next,%
\begin{equation*}
-\left( a_{i,j}v_{x_{i}x_{j}}v_{t}+a_{i,j}v_{x_{j}x_{i}}v_{t}\right)
=-\left( a_{i,j}v_{x_{i}}v_{t}\right)
_{x_{j}}+a_{i,j}v_{x_{i}}v_{tx_{j}}+\left( a_{i,j}\right)
_{x_{j}}v_{x_{i}}v_{t}
\end{equation*}%
\begin{equation*}
-\left( a_{i,j}v_{x_{j}}v_{t}\right)
_{x_{i}}+a_{i,j}v_{tx_{i}}v_{x_{j}}+\left( a_{i,j}\right)
_{x_{i}}v_{x_{j}}v_{t}
\end{equation*}%
\begin{equation*}
=\left( a_{i,j}v_{x_{i}}v_{x_{j}}\right) _{t}-\left( a_{i,j}\right)
_{t}v_{x_{i}}v_{x_{j}}+\left( a_{i,j}\right) _{x_{j}}v_{x_{i}}v_{t}+\left(
a_{i,j}\right) _{x_{i}}v_{x_{j}}v_{t}
\end{equation*}%
\begin{equation}
-\left[ \left( a_{i,j}v_{x_{i}}v_{t}\right) _{x_{j}}+\left(
a_{i,j}v_{x_{j}}v_{t}\right) _{x_{i}}\right]  \label{2.3}
\end{equation}%
\begin{equation*}
=\left( a_{i,j}v_{x_{i}}v_{x_{j}}\right) _{t}-\left( a_{i,j}\right)
_{t}v_{x_{i}}v_{x_{j}}+\left( a_{i,j}\right) _{x_{j}}v_{x_{i}}v_{t}+\left(
a_{i,j}\right) _{x_{i}}v_{x_{j}}v_{t}
\end{equation*}%
\begin{equation*}
-\left[ \left( a_{i,j}v_{x_{i}}v_{t}\right) _{x_{j}}+\left(
a_{i,j}v_{x_{j}}v_{t}\right) _{x_{i}}\right] .
\end{equation*}%
Applying the Cauchy-Schwarz inequality \textquotedblleft with $\varepsilon "$
to (\ref{2.3}),%
\begin{equation}
2ab\geq -\varepsilon a^{2}-\frac{1}{\varepsilon }b^{2},\forall a,b\in 
\mathbb{R},\forall \varepsilon >0,  \label{2.30}
\end{equation}%
and using (\ref{2.2}), we obtain%
\begin{equation*}
-2v_{t}Lv\geq -C\left( \nabla v\right) ^{2}-\frac{1}{2}v_{t}^{2}+\dsum%
\limits_{i,j=1}^{n}\left( a_{i,j}v_{x_{i}}v_{x_{j}}\right) _{t}
\end{equation*}%
\begin{equation*}
+\dsum\limits_{i,j=1}^{n}\left[ \left( a_{i,j}v_{x_{i}}v_{t}\right)
_{x_{j}}+\left( a_{i,j}v_{x_{j}}v_{t}\right) _{x_{i}}\right] .
\end{equation*}%
Hence,%
\begin{equation*}
v_{t}^{2}-2v_{t}Lv\geq \frac{1}{2}v_{t}^{2}-C\left( \nabla v\right) ^{2}
\end{equation*}%
\begin{equation}
+\left( \dsum\limits_{i,j=1}^{n}a_{i,j}v_{x_{i}}v_{x_{j}}\right)
_{t}-\dsum\limits_{i,j=1}^{n}\left[ \left( a_{i,j}v_{x_{i}}v_{t}\right)
_{x_{j}}+\left( a_{i,j}v_{x_{j}}v_{t}\right) _{x_{i}}\right] .  \label{2.4}
\end{equation}

\textbf{Step 2}. Estimate the term $-2\nu \left( t+1\right) ^{\nu -1}v_{t}v$
in the third line of (\ref{2.1}). We have:%
\begin{equation}
-2\nu \left( t+1\right) ^{\nu -1}v_{t}v=\left( -\nu \left( t+1\right) ^{\nu
-1}v^{2}\right) _{t}+\nu \left( \nu -1\right) \left( t+1\right) ^{\nu
-2}v^{2}  \label{2.5}
\end{equation}%
\begin{equation*}
\geq C\nu ^{2}\left( t+1\right) ^{\nu -1}v^{2}+\left( -\nu \left( t+1\right)
^{\nu -1}v^{2}\right) _{t}.
\end{equation*}

\textbf{Step 3}. Sum up (\ref{2.4}), (\ref{2.5})$.$ Then, taking into
account (\ref{2.1}), we obtain%
\begin{equation*}
\left( u_{t}-Lu\right) ^{2}\exp \left( 2\left( t+1\right) ^{\nu }\right)
\geq \frac{1}{2}v_{t}^{2}-C\left( \nabla v\right) ^{2}+C\nu ^{2}\left(
t+1\right) ^{\nu -1}v^{2}
\end{equation*}%
\begin{equation*}
+\left( \dsum\limits_{i,j=1}^{n}a_{i,j}v_{x_{i}}v_{x_{j}}-\nu \left(
t+1\right) ^{\nu -1}v^{2}\right) _{t}+\dsum\limits_{i,j=1}^{n}\left[ \left(
-a_{i,j}v_{x_{i}}v_{t}\right) _{x_{j}}+\left( -a_{i,j}v_{x_{j}}v_{t}\right)
_{x_{i}}\right] .
\end{equation*}%
Replacing here $v$ with $u=v\exp \left( -\left( t+1\right) ^{\nu }\right) $
and using $v_{t}^{2}/2\geq 0,$ we obtain%
\begin{equation*}
\left( u_{t}-Lu\right) ^{2}\exp \left( 2\lambda \left( t+1\right) ^{\nu
}\right) \geq C\left[ -\left( \nabla u\right) ^{2}+\nu ^{2}\left( t+1\right)
^{\nu -1}u^{2}\right] \exp \left( 2\left( t+1\right) ^{\nu }\right)
\end{equation*}%
\begin{equation}
+\left[ \left( \dsum\limits_{i,j=1}^{n}a_{i,j}u_{x_{i}}u_{x_{j}}-\nu \left(
t+1\right) ^{\nu -1}u^{2}\right) \exp \left( 2\left( t+1\right) ^{\nu
}\right) \right] _{t}  \label{2.6}
\end{equation}%
\begin{equation*}
+\dsum\limits_{i,j=1}^{n}\left\{ \left[ -a_{i,j}u_{x_{i}}\left( u_{t}+\nu
\left( t+1\right) ^{\nu -1}u\right) \exp \left( 2\left( t+1\right) ^{\nu
}\right) \right] _{x_{j}}\right\}
\end{equation*}%
\begin{equation*}
+\dsum\limits_{i,j=1}^{n}\left\{ \left[ -a_{i,j}u_{x_{j}}\left( u_{t}+\nu
\left( t+1\right) ^{\nu -1}u\right) \exp \left( 2\left( t+1\right) ^{\nu
}\right) \right] _{x_{i}}\right\} .
\end{equation*}%
What is not good in estimate (\ref{2.6}) is that the first line in its right
hand side contains both positive and negative terms. However, only positive
terms must be in such cases in any Carleman estimate. Hence, we continue
with further steps.

\textbf{Step 4}. Estimate from the below the expression $\left(
u_{t}-Lu\right) u\exp \left( 2\left( t+1\right) ^{\nu }\right) .$ We have 
\begin{equation*}
\left( u_{t}-Lu\right) u\exp \left( 2\left( t+1\right) ^{\nu }\right)
=\left( \frac{u^{2}}{2}\exp \left( 2\left( t+1\right) ^{\nu }\right) \right)
_{t}-\nu \left( t+1\right) ^{\nu -1}u^{2}\exp \left( 2\left( t+1\right)
^{\nu }\right)
\end{equation*}%
\begin{equation*}
-\frac{1}{2}\dsum\limits_{i,j=1}^{n}\left( a_{ij}u_{x_{i}}u\exp \left(
2\left( t+1\right) ^{\nu }\right) \right) _{x_{j}}-\frac{1}{2}%
\dsum\limits_{i,j=1}^{n}\left( a_{ij}u_{x_{j}}u\exp \left( 2\left(
t+1\right) ^{\nu }\right) \right) _{x_{i}}
\end{equation*}%
\begin{equation}
+\dsum\limits_{i,j=1}^{n}a_{ij}u_{x_{i}}u_{j}\exp \left( 2\left( t+1\right)
^{\nu }\right)  \label{2.7}
\end{equation}%
\begin{equation*}
+\frac{1}{2}\dsum\limits_{i,j=1}^{n}\left( a_{ij}\right)
_{x_{j}}u_{x_{i}}u\exp \left( 2\left( t+1\right) ^{\nu }\right) +\frac{1}{2}%
\dsum\limits_{i,j=1}^{n}\left( a_{ij}\right) _{x_{i}}u_{x_{j}}u\exp \left(
2\left( t+1\right) ^{\nu }\right) .
\end{equation*}%
By (\ref{1.3}) 
\begin{equation*}
\dsum\limits_{i,j=1}^{n}a_{ij}u_{x_{i}}u_{j}\exp \left( 2\left( t+1\right)
^{\nu }\right) \geq \mu _{1}\left( \nabla u\right) ^{2}\exp \left( 2\left(
t+1\right) ^{\nu }\right) .
\end{equation*}%
Hence, using (\ref{2.30}) and (\ref{2.7}), we obtain for all $\nu \geq \nu
_{0}$ 
\begin{equation*}
\left( u_{t}-Lu\right) u\exp \left( 2\left( t+1\right) ^{\nu }\right) \geq C 
\left[ \left( \nabla u\right) ^{2}-\nu \left( t+1\right) ^{\nu -1}u^{2}%
\right] \exp \left( 2\left( t+1\right) ^{\nu }\right)
\end{equation*}%
\begin{equation}
+\left( \frac{u^{2}}{2}\exp \left( 2\left( t+1\right) ^{\nu }\right) \right)
_{t}  \label{2.8}
\end{equation}%
\begin{equation*}
-\frac{1}{2}\dsum\limits_{i,j=1}^{n}\left( a_{ij}u_{x_{i}}u\exp \left(
2\left( t+1\right) ^{\nu }\right) \right) _{x_{j}}-\frac{1}{2}%
\dsum\limits_{i,j=1}^{n}\left( a_{ij}u_{x_{j}}u\exp \left( 2\left(
t+1\right) ^{\nu }\right) \right) _{x_{i}}.
\end{equation*}

\textbf{Step 5}. Multiply (\ref{2.8}) by $\sqrt{\nu }$ and sum up with (\ref%
{2.6}). Since $\sqrt{\nu }>>1$ and $\nu ^{2}>>\nu ^{3/2}$ for sufficiently
large $\nu ,$ then we obtain for $\nu \geq \nu _{0}$ 
\begin{equation*}
\left( u_{t}-Lu\right) ^{2}\exp \left( 2\left( t+1\right) ^{\nu }\right) +%
\sqrt{\nu }\left( u_{t}-Lu\right) u\exp \left( 2\left( t+1\right) ^{\nu
}\right)
\end{equation*}%
\begin{equation*}
\geq C\left( \sqrt{\nu }\left( \nabla u\right) ^{2}+\nu ^{2}\left(
t+1\right) ^{\nu -1}u^{2}\right) \exp \left( 2\left( t+1\right) ^{\nu
}\right)
\end{equation*}%
\begin{equation}
\left[ \left( \dsum\limits_{i,j=1}^{n}a_{i,j}u_{x_{i}}u_{x_{j}}-\nu \left(
t+1\right) ^{\nu -1}u^{2}+\sqrt{\nu }\frac{u^{2}}{2}\right) \exp \left(
2\left( t+1\right) ^{\nu }\right) \right] _{t}  \label{2.9}
\end{equation}%
\begin{equation*}
+\dsum\limits_{i,j=1}^{n}\left\{ \left[ a_{i,j}u_{x_{i}}\left( u_{t}+\nu
\left( t+1\right) ^{\nu -1}u\right) \exp \left( 2\left( t+1\right) ^{\nu
}\right) \right] _{x_{j}}\right\}
\end{equation*}%
\begin{equation*}
+\dsum\limits_{i,j=1}^{n}\left\{ \left[ a_{i,j}u_{x_{j}}\left( u_{t}+\nu
\left( t+1\right) ^{\nu -1}u\right) \exp \left( 2\left( t+1\right) ^{\nu
}\right) \right] _{x_{i}}\right\}
\end{equation*}%
\begin{equation*}
-\frac{\sqrt{\nu }}{2}\dsum\limits_{i,j=1}^{n}\left( a_{ij}u_{x_{i}}u\exp
\left( 2\left( t+1\right) ^{\nu }\right) \right) _{x_{j}}-\frac{\sqrt{\nu }}{%
2}\dsum\limits_{i,j=1}^{n}\left( a_{ij}u_{x_{j}}u\exp \left( 2\left(
t+1\right) ^{\nu }\right) \right) _{x_{i}}.
\end{equation*}%
Next, we estimate from the above the left hand side of inequality (\ref{2.9}%
) as 
\begin{equation*}
\left( u_{t}-Lu\right) ^{2}\exp \left( 2\left( t+1\right) ^{\nu }\right) +%
\sqrt{\nu }\left( u_{t}-Lu\right) u\exp \left( 2\left( t+1\right) ^{\nu
}\right)
\end{equation*}%
\begin{equation*}
\leq \frac{3}{2}\left( u_{t}-Lu\right) ^{2}\exp \left( 2\left( t+1\right)
^{\nu }\right) +\frac{\nu }{2}u^{2}\exp \left( 2\left( t+1\right) ^{\nu
}\right) .
\end{equation*}%
Comparing this with (\ref{2.9}), we obtain%
\begin{equation*}
\left( u_{t}-Lu\right) ^{2}\exp \left( 2\left( t+1\right) ^{\nu }\right)
\end{equation*}%
\begin{equation}
\geq C\left( \sqrt{\nu }\left( \nabla u\right) ^{2}+\nu ^{2}\left(
t+1\right) ^{\nu -1}u^{2}\right) \exp \left( 2\left( t+1\right) ^{\nu
}\right)  \label{2.10}
\end{equation}%
\begin{equation*}
+\left[ \left( \dsum\limits_{i,j=1}^{n}a_{i,j}u_{x_{i}}u_{x_{j}}-\nu \left(
t+1\right) ^{\nu -1}u^{2}+\sqrt{\nu }\frac{u^{2}}{2}\right) \exp \left(
2\left( t+1\right) ^{\nu }\right) \right] _{t}
\end{equation*}%
\begin{equation*}
+\func{div}U,\text{ }\nu \geq \nu _{0},
\end{equation*}%
where the vector function $U$ is such that 
\begin{equation}
U\mid _{S_{T}}=0.  \label{2.11}
\end{equation}%
Condition (\ref{2.11}) follows from the boundary condition in (\ref{1}) and
the lines 5 and 6 of (\ref{2.9}).

\textbf{Step 6}. Integrate the pointwise Carleman estimate (\ref{2.10}) over
the domain $Q_{T}.$ Using Gauss' formula and (\ref{2.11}), we obtain (\ref%
{1.9}). $\square $

\section{ Stability Estimate}

\label{sec:3}

For an arbitrary $\tau \in \left( 0,T\right) ,$ denote 
\begin{equation*}
H^{1,0}\left( Q_{T\tau }\right) =\left\{ u:\left\Vert u\right\Vert
_{H^{1,0}\left( Q_{T\tau }\right) }=\left[ \dint\limits_{Q_{T\tau }}\left(
\left( \nabla u\right) ^{2}+u^{2}\right) \left( x,t\right) dxdt\right]
^{1/2}<\infty \right\} .
\end{equation*}

Prior establishing our stability estimate, we prove Lemma 3.1.

\textbf{Lemma 3.1.} \emph{Let }$\delta \in \left( 0,1\right) $\emph{\ be a
sufficiently small number and let the number }$k>0$\emph{. Choose a
sufficiently large number }$\nu =\nu \left( \delta \right) $\emph{\ such
that }%
\begin{equation}
\exp \left( k\left( T+1\right) ^{\nu \left( \delta \right) }\right) =\frac{1%
}{\delta }.  \label{2}
\end{equation}%
\emph{Then for any }$\tau \in \left( 0,T\right) $\emph{\ and for any number }%
$y>0$\emph{\ } 
\begin{equation}
\lim_{\delta \rightarrow 0}\frac{\exp \left( -2\left( \tau +1\right) ^{\nu
\left( \delta \right) }\right) }{\delta ^{y}}=\infty ,  \label{3}
\end{equation}%
\begin{equation}
\lim_{\delta \rightarrow 0}\frac{\exp \left( -2\left( \tau +1\right) ^{\nu
\left( \delta \right) }\right) }{\left( \ln \left( \delta ^{-1}\right)
\right) ^{-y}}=0.  \label{4}
\end{equation}

\textbf{Remark 3.1}. \emph{It follows (\ref{3}) and (\ref{4}) that any
stability estimate via }$\exp \left( -2\left( \tau +1\right) ^{\nu }\right) $%
\emph{\ is weaker than H\"{o}lder and stronger than logarithmic stability
estimate.}

\textbf{Proof of Lemma 3.1.} By (\ref{2})\emph{\ }$k\left( T+1\right) ^{\nu
\left( \delta \right) }=\ln \left( \delta ^{-1}\right) .$ Hence, $\left(
T+1\right) ^{\nu \left( \delta \right) }=\ln \left( \delta ^{-1/k}\right) .$
Hence, 
\begin{equation}
\nu \left( \delta \right) =\ln \left[ \ln \left( \delta ^{-1/k}\right)
^{1/\ln \left( T+1\right) }\right] .  \label{5}
\end{equation}%
Next,%
\begin{equation}
\left( \tau +1\right) ^{\nu \left( \delta \right) }=\exp \left( \ln \left(
\tau +1\right) ^{\nu \left( \delta \right) }\right) =\exp \left( \nu \left(
\delta \right) \ln \left( \tau +1\right) \right) .  \label{6}
\end{equation}%
By (\ref{5})%
\begin{equation}
\nu \left( \delta \right) \ln \left( \tau +1\right) =\ln \left[ \ln \left(
\delta ^{-1/k}\right) ^{\ln \left( \tau +1\right) /\ln \left( T+1\right) }%
\right] =\ln \left[ \ln \left( \delta ^{-1/k}\right) ^{c}\right] ,  \label{7}
\end{equation}%
\begin{equation}
c=c\left( \tau ,T\right) =\frac{\ln \left( \tau +1\right) }{\ln \left(
T+1\right) }\in \left( 0,1\right) .  \label{8}
\end{equation}%
Using (\ref{6})-(\ref{8}), we obtain%
\begin{equation*}
\left( \tau +1\right) ^{\nu \left( \delta \right) }=\ln \left( \delta
^{-1/k}\right) ^{c}=\frac{1}{k^{c}}\ln \left( \delta ^{-1}\right) ^{c}.
\end{equation*}%
Hence,%
\begin{equation}
\exp \left( -2\left( \tau +1\right) ^{\nu \left( \delta \right) }\right)
=\exp \left[ -\frac{2}{k^{c}}\ln \left( \delta ^{-1}\right) ^{c}\right] .
\label{9}
\end{equation}%
Using (\ref{9}), we now prove (\ref{3}). Indeed, 
\begin{equation}
\frac{\exp \left[ -2\ln \left( \delta ^{-1}\right) ^{c}/k^{c}\right] }{%
\delta ^{y}}=\frac{\exp \left[ -2\ln \left( \delta ^{-1}\right) ^{c}/k^{c}%
\right] }{\exp \left( -y\ln \left( \delta ^{-1}\right) \right) }=  \label{90}
\end{equation}%
\begin{equation*}
=\exp \left[ y\ln \left( \delta ^{-1}\right) \left[ 1-\frac{2}{yk^{c}}\ln
\left( \delta ^{-1}\right) ^{c-1}\right] \right] .
\end{equation*}%
Since by (\ref{8}) $c\in \left( 0,1\right) ,$ then 
\begin{equation*}
\lim_{\delta \rightarrow 0}\frac{2}{yk^{c}}\ln \left( \delta ^{-1}\right)
^{c-1}=0.
\end{equation*}%
Hence, 
\begin{equation*}
\lim_{\delta \rightarrow 0}\left\{ \exp \left[ y\ln \left( \delta
^{-1}\right) \left[ 1-\frac{2}{yk^{c}}\ln \left( \delta ^{-1}\right) ^{c-1}%
\right] \right] \right\} =\infty .
\end{equation*}%
This and (\ref{90}) prove (\ref{3}).

We now prove (\ref{4}), which is equivalent with%
\begin{equation}
\lim_{\delta \rightarrow 0}\frac{\exp \left[ -2\ln \left( \delta
^{-1}\right) ^{c}/k^{c}\right] }{\exp \left[ -y\ln \left( \ln \left( \delta
^{-1}\right) \right) \right] }=0.  \label{10}
\end{equation}%
Next,%
\begin{equation*}
\frac{\exp \left[ -2\ln \left( \delta ^{-1}\right) ^{c}/k^{c}\right] }{\exp %
\left[ -y\ln \left( \ln \left( \delta ^{-1}\right) \right) \right] }=\exp %
\left[ -\frac{2}{k^{c}}\ln \left( \delta ^{-1}\right) ^{c}+y\ln \left( \ln
\left( \delta ^{-1}\right) \right) \right] ,
\end{equation*}%
\begin{equation}
-\frac{2}{k^{c}}\ln \left( \delta ^{-1}\right) ^{c}+y\ln \left( \ln \left(
\delta ^{-1}\right) \right) =-\frac{2}{k^{c}}\ln \left( \delta ^{-1}\right)
^{c}\left[ 1-\frac{yk^{c}}{2}\frac{\ln \left( \ln \left( \delta ^{-1}\right)
\right) }{\ln \left( \delta ^{-1}\right) ^{c}}\right] .  \label{11}
\end{equation}%
Obviously%
\begin{equation}
\lim_{\delta \rightarrow 0}\frac{\ln \left( \ln \left( \delta ^{-1}\right)
\right) }{\ln \left( \delta ^{-1}\right) ^{c}}=0.  \label{12}
\end{equation}%
Thus, (\ref{10}) follows from (\ref{11}) and (\ref{12}). $\square $

\textbf{Theorem 3.1} (stability estimate). \emph{Assume that conditions (\ref%
{1.1})-(\ref{1.5}) are in place. \ Suppose that two functions }$%
u_{1},u_{2}\in H^{2}\left( Q_{T}\right) $\emph{\ are solutions of problem (%
\ref{1.6}), (\ref{1.07}) with different data at }$\left\{ t=T\right\} $\emph{%
, }%
\begin{equation}
u_{1}\left( x,T\right) =g_{1}\left( x\right) ,\text{ }u_{2}\left( x,T\right)
=g_{2}\left( x\right) ,\text{ }f\left( x\right) =g_{1}\left( x\right)
-g_{2}\left( x\right) .  \label{14}
\end{equation}%
\emph{Suppose that }%
\begin{equation}
\left\Vert f\right\Vert _{L_{2}\left( \Omega \right) }\leq \delta ,
\label{2.12}
\end{equation}%
\emph{\ }%
\begin{equation}
\left\Vert \nabla u_{i}\left( x,0\right) \right\Vert _{L_{2}\left( \Omega
\right) }\leq M,i=1,2,  \label{2.13}
\end{equation}%
\emph{where }$M>0$\emph{\ is a number and the parameter }$\delta $\emph{\
characterizes the level of noise in the data }$u\left( x,T\right) $\emph{.
Denote }$w=u_{1}-u_{2}.$ \emph{Let }$\overline{C}>0$\emph{\ be the number in
(\ref{1.10}) and }$\nu _{0}$\emph{\ be the parameter of Theorem 2.1. Then
there exist constants }%
\begin{equation*}
C_{1}=C_{1}\left( \mu _{1},\mu _{2},\max_{i,j}\left\Vert a_{ij}\right\Vert
_{C^{1}\left( \overline{Q}_{T}\right) },Q_{T},\overline{C}\right) >0,
\end{equation*}%
\begin{equation*}
\nu _{1}=\nu _{1}\left( \mu _{1},\mu _{2},\max_{i,j}\left\Vert
a_{ij}\right\Vert _{C^{1}\left( \overline{Q}_{T}\right) },Q_{T},\overline{C}%
\right) \geq \nu _{0}
\end{equation*}%
\emph{\ depending only on listed parameters such that if the number }$\delta
_{0}$\emph{\ is so small that} 
\begin{equation}
\ln \left[ \ln \left( \delta _{0}^{-1/3}\right) ^{1/\ln \left( T+1\right) }%
\right] \geq \nu _{1},  \label{2.14}
\end{equation}%
\emph{then the following estimate holds for any }$\tau \in \left( 0,T\right) 
$\emph{\ and for all }$\delta \in \left( 0,\delta _{0}\right) $ \emph{\ }%
\begin{equation}
\left\Vert w\right\Vert _{H^{1,0}\left( Q_{T\tau }\right) }\leq C_{1}\left(
M+1\right) \exp \left[ -\frac{1}{3^{c}}\ln \left( \delta ^{-1}\right) ^{c}%
\right] ,  \label{30.2}
\end{equation}%
\emph{where} \emph{the constant }$C_{1}$\emph{\ is independent on }$M$\emph{%
\ and the number }$c=c\left( \tau ,T\right) \in \left( 0,1\right) $ \emph{is
defined in (\ref{8}).}

Estimate (\ref{30.2}) is between H\"{o}lder and logarithmic stability
estimates, see Lemma 3.1 and Remark 3.1. Everywhere below $C_{1}>0$ denotes
different numbers depending on the same parameters as ones listed in this
theorem.

\textbf{Proof}. It follows from (\ref{1.10})-(\ref{1.8}) and (\ref{14}) that%
\begin{equation}
\left\vert w_{t}-Lw\right\vert \leq \overline{C}\left( \left\vert \nabla
w\right\vert +\left\vert w\right\vert \right) \text{ a.e. in }Q_{T},
\label{3.1}
\end{equation}%
\begin{equation}
w\mid _{S_{T}}=0,  \label{3.2}
\end{equation}%
\begin{equation}
w\left( x,T\right) =f\left( x\right) .  \label{3.3}
\end{equation}%
Square both sides of inequality (\ref{3.1}), then multiply by $\exp \left(
2\left( t+1\right) ^{\nu }\right) $ integrate over $Q_{T}$ and then apply
Theorem 2.1 taking into account (\ref{3.2}) and (\ref{3.3}). We obtain%
\begin{equation*}
2\overline{C}^{2}\dint\limits_{Q_{T}}\left( \left( \nabla w\right)
^{2}+w^{2}\right) \exp \left( 2\left( t+1\right) ^{\nu }\right) dxdt
\end{equation*}%
\begin{equation*}
\geq \dint\limits_{Q_{T}}\left( w_{t}-Lw\right) ^{2}\exp \left( 2\left(
t+1\right) ^{\nu }\right) dxdt
\end{equation*}%
\begin{equation}
\geq C\sqrt{\nu }\dint\limits_{Q_{T}}\left( \nabla w\right) ^{2}\exp \left(
2\left( t+1\right) ^{\nu }\right) dxdt+C\nu ^{2}\dint\limits_{Q_{T\tau
}}w^{2}\exp \left( 2\left( t+1\right) ^{\nu }\right) dxdt  \label{3.4}
\end{equation}%
\begin{equation*}
-C\exp \left( 3\left( T+1\right) ^{\nu }\right) \left\Vert f\right\Vert
_{L_{2}\left( \Omega \right) }^{2}-C\left\Vert \nabla w\left( x,0\right)
\right\Vert _{L_{2}\left( \Omega \right) }^{2},\forall \nu \geq \nu _{0}.
\end{equation*}%
Choose $\nu _{1}=\nu _{1}\left( \mu _{1},\mu _{2},\max_{i,j}\left\Vert
a_{ij}\right\Vert _{C^{1}\left( \overline{Q}_{T}\right) },Q_{T},\overline{C}%
\right) \geq \nu _{0}>1$ so that $C\sqrt{\nu _{1}}/4\geq \overline{C}^{2}.$
Then (\ref{2.12}) and (\ref{3.4}) imply that for all $\nu \geq \nu _{1}$ 
\begin{equation*}
\delta ^{2}\exp \left( 3\left( T+1\right) ^{\nu }\right) +\left\Vert \nabla
w\left( x,0\right) \right\Vert _{L_{2}\left( \Omega \right) }^{2}
\end{equation*}%
\begin{equation}
\geq C_{1}\sqrt{\nu }\dint\limits_{Q_{T}}\left( \nabla w\right) ^{2}\exp
\left( 2\left( t+1\right) ^{\nu }\right) dxdt+C_{1}\nu
^{2}\dint\limits_{Q_{T}}w^{2}\exp \left( 2\left( t+1\right) ^{\nu }\right)
dxdt  \label{3.5}
\end{equation}%
\begin{equation*}
\geq C_{1}\sqrt{\nu }\dint\limits_{Q_{T\tau }}\left( \nabla w\right)
^{2}\exp \left( 2\left( t+1\right) ^{\nu }\right) dxdt+C_{1}\nu
^{2}\dint\limits_{Q_{T\tau }}w^{2}\exp \left( 2\left( t+1\right) ^{\nu
}\right) dxdt
\end{equation*}%
\begin{equation*}
\geq C_{1}\exp \left( 2\left( \tau +1\right) ^{\nu }\right) \left\Vert
w\right\Vert _{H^{1,0}\left( Q_{T\tau }\right) }^{2}.
\end{equation*}%
Hence,%
\begin{equation}
\left\Vert w\right\Vert _{H^{1,0}\left( Q_{T\tau }\right) }^{2}\leq
C_{1}\delta ^{2}\exp \left( 3\left( T+1\right) ^{\nu }\right)
+C_{1}\left\Vert \nabla w\left( x,0\right) \right\Vert _{L_{2}\left( \Omega
\right) }^{2}\exp \left( -2\left( \tau +1\right) ^{\nu }\right) .
\label{3.6}
\end{equation}%
Choose $\nu =\nu \left( \delta \right) $ such that (\ref{2}) would be
satisfied with $k=3$, i.e. 
\begin{equation}
\exp \left( 3\left( T+1\right) ^{\nu }\right) =\frac{1}{\delta }.
\label{3.7}
\end{equation}%
Hence, 
\begin{equation}
\nu =\nu \left( \delta \right) =\ln \left[ \left( \ln \left( \delta
^{-1/3}\right) ^{1/\ln \left( T+1\right) }\right) \right] .  \label{3.8}
\end{equation}%
The choice (\ref{3.8}) is possible since (\ref{2.14}) holds and $\delta \in
\left( 0,\delta _{0}\right) .$ Hence, by (\ref{9}) 
\begin{equation}
\exp \left( -2\left( \tau +1\right) ^{\nu \left( \delta \right) }\right)
=\exp \left[ -\frac{2}{3^{c}}\ln \left( \delta ^{-1}\right) ^{c}\right] ,
\label{3.9}
\end{equation}
Hence, (\ref{2.13}), (\ref{3.6}), (\ref{3.7}) and (\ref{3.9}) imply that 
\begin{equation}
\left\Vert w\right\Vert _{H^{1,0}\left( Q_{T\tau }\right) }^{2}\leq
C_{1}\delta +C_{1}M^{2}\exp \left[ -\frac{2}{3^{c}}\ln \left( \delta
^{-1}\right) ^{c}\right] .  \label{3.10}
\end{equation}%
The target estimate (\ref{30.2}) of this theorem obviously follows from (\ref%
{3}) and (\ref{3.10}). $\square $

\section{The Quasi-Reversibility Method for the Linear Case}

\label{sec:4}

We assume in this section that the function $F\left( \nabla u,u,x,t\right) $
in (\ref{1.6}) is linear with respect to the function $u$ and its first
derivatives,%
\begin{equation}
F\left( \nabla u,u,x,t\right) =Au+p\left( x,t\right)
=\dsum\limits_{j=1}^{n}b_{j}\left( x,t\right) u_{x_{j}}+c\left( x,t\right)
u+p\left( x,t\right) ,  \label{4.1}
\end{equation}%
where functions $b_{j},c,p\in C\left( \overline{Q}_{T}\right) .$ Then (\ref%
{1.6})-(\ref{1.8}) become 
\begin{equation}
u_{t}=Lu+Au+p\left( x,t\right) ,\left( x,t\right) \in Q_{T},  \label{4.2}
\end{equation}%
\begin{equation}
u\mid _{S_{T}}=0,  \label{4.4}
\end{equation}%
\begin{equation}
u\left( x,T\right) =g\left( x\right) ,x\in \Omega .  \label{4.5}
\end{equation}%
Assuming that $g\left( x\right) =0$ for $x\in \partial \Omega $ and 
\begin{equation}
g\in H^{2}\left( \Omega \right) ,  \label{4.6}
\end{equation}%
consider the function $v\left( x,t\right) =u\left( x,t\right) -g\left(
x\right) .$ Then (\ref{4.2})-(\ref{4.6}) lead to%
\begin{equation}
v_{t}=Lv+Av+q\left( x,t\right) ,\left( x,t\right) \in Q_{T},  \label{4.7}
\end{equation}%
\begin{equation}
v\mid _{S_{T}}=0,  \label{4.8}
\end{equation}%
\begin{equation}
v\left( x,T\right) =0,x\in \Omega ,  \label{4.9}
\end{equation}%
\begin{equation}
q\left( x,t\right) =Lg+Ag+p\left( x,t\right) \in L_{2}\left( Q_{T}\right) .
\label{4.10}
\end{equation}

We introduce the subspace $H_{0}^{2}\left( Q_{T}\right) $ of the space $%
H^{2}\left( Q_{T}\right) $ as%
\begin{equation*}
H_{0}^{2}\left( Q_{T}\right) =\left\{ w\in H^{2}\left( Q_{T}\right) :w\mid
_{S_{T}}=0,w\left( x,T\right) =0\right\} .
\end{equation*}

The QRM for problem (\ref{4.7})-(\ref{4.10}) amounts to the minimization of
the following functional $J_{\alpha },$%
\begin{equation}
J_{\alpha }\left( v\right) =\dint\limits_{Q_{T}}\left( v_{t}-Lv-Av-q\right)
^{2}dxdt+\alpha \left\Vert v\right\Vert _{H^{2}\left( Q_{T}\right) }^{2},
\label{4.11}
\end{equation}%
where $\alpha \in \left( 0,1\right) $ is the regularization parameter. We
arrive at the following problem:

\textbf{Minimization Problem 1}. \emph{Minimize the functional }$J_{\alpha
}\left( v\right) $\emph{\ on the space} $H_{0}^{2}\left( Q_{T}\right) .$

\textbf{Theorem 4.1}. \emph{Assume that conditions (\ref{1.1})-(\ref{1.5}), (%
\ref{4.10}) hold. Then there exists unique minimizer }$v_{\min }\in
H_{0}^{2}\left( Q_{T}\right) $\emph{\ of the functional} $J_{\alpha }\left(
v\right) .$

\textbf{Proof}. Let $\left[ .,.\right] $ denotes the scalar product in $%
H^{2}\left( Q_{T}\right) .$ By the variational principle, any minimizer $%
v_{\min }\in H_{0}^{2}\left( Q_{T}\right) ,$ if it exists, satisfies the
following integral identity%
\begin{equation}
\dint\limits_{Q_{T}}\left( \partial _{t}v_{\min }-Lv_{\min }-Av_{\min
}\right) \left( h_{t}-Lh-Ah\right) dxdt+\alpha \left[ v_{\min },h\right]
\label{4.12}
\end{equation}%
\begin{equation*}
=\dint\limits_{Q_{T}}q\left( h_{t}-Lh-Ah\right) dxdt,\forall h\in
H_{0}^{2}\left( Q_{T}\right) .
\end{equation*}%
Since 
\begin{equation*}
\dint\limits_{Q_{T}}\left( v_{t}-Lv-Av\right) ^{2}dxdt+\alpha \left[ v,v%
\right] ^{2}\geq \alpha \left[ v,v\right] ^{2}=\alpha \left\Vert
v\right\Vert _{H^{2}\left( Q_{T}\right) }^{2},\text{ }\forall v\in
H_{0}^{2}\left( Q_{T}\right) ,
\end{equation*}%
\begin{equation*}
\dint\limits_{Q_{T}}\left( v_{t}-Lv-Av\right) ^{2}dxdt+\alpha \left[ v,v%
\right] ^{2}\leq C_{1}\left\Vert v\right\Vert _{H^{2}\left( Q_{T}\right)
}^{2},\text{ }\forall v\in H_{0}^{2}\left( Q_{T}\right) ,
\end{equation*}%
then the equality%
\begin{equation*}
\left\{ v,h\right\} =\dint\limits_{Q_{T}}\left( v_{t}-Lv-Mv\right) \left(
h_{t}-Lh-Mh\right) dxdt+\alpha \left[ v,h\right] ,\forall v,h\in
H_{0}^{2}\left( Q_{T}\right) .
\end{equation*}%
defines a new scalar product in the space $H_{0}^{2}\left( Q_{T}\right) .$
We rewrite integral identity (\ref{4.12}) as%
\begin{equation}
\left\{ v_{\min },h\right\} =\dint\limits_{Q_{T}}q\left( h_{t}-Lh-Ah\right)
dxdt.  \label{4.13}
\end{equation}%
Next, by the Cauchy-Schwarz inequality%
\begin{equation*}
\dint\limits_{Q_{T}}q\left( h_{t}-Lh-Ah\right) dxdt\leq \left\Vert
q\right\Vert _{L_{2}\left( Q_{T}\right) }\left\Vert h_{t}-Lh-Ah\right\Vert
_{L_{2}\left( Q_{T}\right) }\leq C_{1}\left\Vert q\right\Vert _{L_{2}\left(
Q_{T}\right) }\left\Vert h\right\Vert _{H^{2}\left( Q_{T}\right) }.
\end{equation*}%
Hence, by Riesz theorem, there exists a unique function $p\in
H_{0}^{2}\left( Q_{T}\right) $ such that 
\begin{equation*}
\dint\limits_{Q_{T}}q\left( h_{t}-Lh-Ah\right) dxdt=\left\{ p,h\right\} ,%
\text{ }\forall h\in H_{0}^{2}\left( Q_{T}\right) .
\end{equation*}%
Comparing this with (\ref{4.13}), we obtain%
\begin{equation*}
\left\{ v_{\min },h\right\} =\left\{ p,h\right\} ,\text{ }\forall h\in
H_{0}^{2}\left( Q_{T}\right) .
\end{equation*}%
Therefore, there exists unique minimizer $v_{\min }=p\in H_{0}^{2}\left(
Q_{T}\right) $ of the functional $J_{\alpha }\left( v\right) .$ $\square $

In the regularization theory, the function $v_{\min }$ is called the \emph{%
regularized solution} of problem (\ref{4.7})-(\ref{4.10}) \cite{BK,T}. The
next step after Theorem 4.1 is to prove convergence of regularized solutions
to the exact one when the noise in the data tends to zero. While we have
used only Riesz theorem to prove existence and uniqueness of the minimizer,
the convergence result requires the Carleman estimate of Theorem 2.1.

Let the function $q^{\ast }\in L_{2}\left( Q_{T}\right) $ be the exact data
in problem (\ref{4.7})-(\ref{4.10}), i.e. the data without a noise in it.
Suppose that there exists the exact solution $v^{\ast }\in H_{0}^{2}\left(
Q_{T}\right) $ of problem (\ref{4.7})-(\ref{4.10}) with $q=q^{\ast }.$
Theorem 1.2 implies that this solution is unique. Let $q\in L_{2}\left(
Q_{T}\right) $ be the noisy data and let $v_{\min }\in H_{0}^{2}\left(
Q_{T}\right) $ be the corresponding minimizer of functional (\ref{4.11})
(Theorem 4.1). We assume that 
\begin{equation}
\left\Vert q-q^{\ast }\right\Vert _{L_{2}\left( Q_{T}\right) }<\delta ,
\label{4.14}
\end{equation}%
where $\delta >0$ is the noise level. Consider the difference 
\begin{equation}
w_{\delta }=v_{\min }-v^{\ast }.  \label{4.1400}
\end{equation}

Theorem 4.2 estimates the function $w$ via $\delta .$ It follows from the
regularization theory that we need to assume a certain dependence $\alpha
=\alpha \left( \delta \right) $ of the regularization parameter on the noise
level $\delta .$

\textbf{Theorem 4.2} (convergence rate). \emph{Assume that conditions (\ref%
{1.1})-(\ref{1.10}) and (\ref{4.14}) are in place. Let }$\nu _{0}>1$ \emph{%
and} $\nu _{1}\geq \nu _{0}$\emph{\ be the numbers of Theorems 2.1 and 3.1
respectively and let }$\alpha =\alpha \left( \delta \right) =\delta ^{2}.$%
\emph{\ Then there exists a number }%
\begin{equation*}
\nu _{2}=\nu _{2}\left( \mu _{1},\mu _{2},\max_{i,j}\left\Vert
a_{ij}\right\Vert _{C^{1}\left( \overline{Q}_{T}\right) },Q_{T},\overline{C}%
\right) \geq \nu _{1}\geq \nu _{0}>1
\end{equation*}%
\emph{depending only on listed parameters} \emph{such that if }$\delta \in
\left( 0,\delta _{0}\right) $ and 
\begin{equation}
\ln \left[ \left( \ln \left( \delta _{0}^{-1/2}\right) ^{1/\ln \left(
T+1\right) }\right) \right] \geq \nu _{2},  \label{4.140}
\end{equation}%
\emph{then the following convergence estimate of the QRM holds for every }$%
\tau \in \left( 0,T\right) :$\emph{\ }%
\begin{equation}
\left\Vert w_{\delta }\right\Vert _{H^{1,0}\left( Q_{T\tau }\right) }\leq
C_{1}\left( 1+\left\Vert v^{\ast }\right\Vert _{H^{2}\left( Q_{T}\right)
}\right) \exp \left[ -\frac{1}{2^{c}}\ln \left( \delta ^{-1}\right) ^{c}%
\right] ,  \label{4.15}
\end{equation}%
\emph{where the function }$w_{\delta }$ \emph{is defined in (\ref{4.1400})
and the number }$c=c\left( \tau ,T\right) \in \left( 0,1\right) $\emph{\ is
defined in (\ref{8}). }

As to estimate (\ref{4.15}), also see Lemma 3.1 and Remark 3.1.

\textbf{Proof of Theorem 4.2}. The function $v^{\ast }\in H_{0}^{2}\left(
Q_{T}\right) $ satisfies the following integral identity%
\begin{equation}
\dint\limits_{Q_{T}}\left( v_{t}^{\ast }-Lv^{\ast }-Av^{\ast }\right) \left(
h_{t}-Lh-Ah\right) dxdt+\alpha \left[ v^{\ast },h\right]  \label{4.16}
\end{equation}%
\begin{equation*}
=\dint\limits_{Q_{T}}q^{\ast }\left( h_{t}-Lh-Ah\right) dxdt+\alpha \left[
v^{\ast },h\right] ,\forall h\in H_{0}^{2}\left( Q_{T}\right) .
\end{equation*}%
Subtracting (\ref{4.16}) from (\ref{4.12}) and using (\ref{4.1400}), we
obtain%
\begin{equation}
\dint\limits_{Q_{T}}\left( w_{\delta t}-Lw_{\delta }-Aw_{\delta }\right)
\left( h_{t}-Lh-Ah\right) dxdt+\alpha \left[ w_{\delta },h\right]
\label{4.17}
\end{equation}%
\begin{equation*}
=\dint\limits_{Q_{T}}\left( q-q^{\ast }\right) \left( h_{t}-Lh-Ah\right)
dxdt-\alpha \left[ v^{\ast },h\right] ,\forall h\in H_{0}^{2,1}\left(
Q_{T}\right) .
\end{equation*}%
Set in (\ref{4.17}) $h=w_{\delta }.$ Using Cauchy-Schwarz inequality, we
obtain%
\begin{equation*}
\dint\limits_{Q_{T}}\left( w_{\delta t}-Lw_{\delta }-Aw_{\delta }\right)
^{2}dxdt+\alpha \left\Vert w_{\delta }\right\Vert _{H^{2}\left( Q_{T}\right)
}^{2}\leq \dint\limits_{Q_{T}}\left\vert q-q^{\ast }\right\vert \cdot
\left\vert w_{\delta t}-Lw_{\delta }-Aw_{\delta }\right\vert dxdt-\alpha 
\left[ v^{\ast },w_{\delta }\right]
\end{equation*}%
\begin{equation*}
\leq \frac{1}{2}\left\Vert q-q^{\ast }\right\Vert _{L_{2}\left( Q\right)
}^{2}+\frac{1}{2}\dint\limits_{Q_{T}}\left( w_{\delta t}-Lw_{\delta
}-Aw_{\delta }\right) ^{2}dxdt+\frac{\alpha }{2}\left\Vert v^{\ast
}\right\Vert _{H^{2,1}\left( Q_{T}\right) }^{2}+\frac{\alpha }{2}\left\Vert
w_{\delta }\right\Vert _{H^{2}\left( Q_{T}\right) }^{2}.
\end{equation*}%
Hence, 
\begin{equation}
\dint\limits_{Q_{T}}\left( w_{\delta t}-Lw_{\delta }-Aw_{\delta }\right)
^{2}dxdt+\alpha \left\Vert w_{\delta }\right\Vert _{H^{2}\left( Q_{T}\right)
}^{2}\leq \delta ^{2}+\alpha \left\Vert v^{\ast }\right\Vert _{H^{2}\left(
Q_{T}\right) }^{2}.  \label{4.18}
\end{equation}%
Since $\alpha =\alpha \left( \delta \right) =\delta ^{2},$ then (\ref{4.18})
implies that $\left\Vert w_{\delta }\right\Vert _{H^{2}\left( Q_{T}\right)
}^{2}\leq 1+\left\Vert v^{\ast }\right\Vert _{H^{2}\left( Q_{T}\right)
}^{2}. $ Hence, trace theorem leads to 
\begin{equation}
\left\Vert \nabla w_{\delta }\left( x,0\right) \right\Vert _{L_{2}\left(
\Omega \right) }\leq C\left( 1+\left\Vert v^{\ast }\right\Vert _{H^{2}\left(
Q_{T}\right) }\right) .  \label{4.180}
\end{equation}

To proceed further, we apply to (\ref{4.18}) the Carleman estimate of
Theorem 2.1. Let $\nu \geq \nu _{1}.$ We have 
\begin{equation*}
\dint\limits_{Q_{T}}\left( w_{\delta t}-Lw_{\delta }-Aw_{\delta }\right)
^{2}dxdt=\dint\limits_{Q_{T}}\exp \left( 2\left( t+1\right) ^{\nu }\right)
\exp \left( -2\left( t+1\right) ^{\nu }\right) \left( w_{\delta
t}-Lw_{\delta }-Aw_{\delta }\right) ^{2}dxdt
\end{equation*}%
\begin{equation*}
\geq \exp \left( -2\left( T+1\right) ^{\nu }\right)
\dint\limits_{Q_{T}}\left( w_{\delta t}-Lw_{\delta }-Aw_{\delta }\right)
^{2}\exp \left( 2\left( t+1\right) ^{\nu }\right) dxdt.
\end{equation*}%
Hence, using (\ref{4.18}), we obtain%
\begin{equation}
\dint\limits_{Q_{T}}\left( w_{\delta t}-Lw_{\delta }-Aw_{\delta }\right)
^{2}\exp \left( 2\left( t+1\right) ^{\nu }\right) dxdt\leq \exp \left(
2\left( T+1\right) ^{\nu }\right) \delta ^{2}\left( 1+\left\Vert v^{\ast
}\right\Vert _{H^{2}\left( Q_{T}\right) }^{2}\right) .  \label{4.19}
\end{equation}%
We have%
\begin{equation}
\dint\limits_{Q_{T}}\left( w_{\delta t}-Lw_{\delta }-Aw_{\delta }\right)
^{2}\exp \left( 2\left( t+1\right) ^{\nu }\right) dxdt  \label{4.20}
\end{equation}%
\begin{equation*}
\geq \dint\limits_{Q_{T}}\left( w_{\delta t}-Lw_{\delta }\right) ^{2}\exp
\left( 2\left( t+1\right) ^{\nu }\right)
dxdt-C_{1}\dint\limits_{Q_{T}}\left( \left( \nabla w_{\delta }\right)
^{2}+w_{\delta }^{2}\right) ^{2}\exp \left( 2\left( t+1\right) ^{\nu
}\right) dxdt.
\end{equation*}%
Next, since $w_{\delta }\left( x,T\right) =0,$ then by (\ref{3.3}) and (\ref%
{3.4})%
\begin{equation*}
\dint\limits_{Q_{T}}\left( w_{\delta t}-Lw_{\delta }\right) ^{2}\exp \left(
2\left( t+1\right) ^{\nu }\right) dxdt-C_{1}\dint\limits_{Q_{T}}\left(
\left( \nabla w_{\delta }\right) ^{2}+w_{\delta }^{2}\right) ^{2}\exp \left(
2\left( t+1\right) ^{\nu }\right) dxdt
\end{equation*}%
\begin{equation}
\geq C\sqrt{\nu }\dint\limits_{Q_{T}}\left( \nabla w_{\delta }\right)
^{2}\exp \left( 2\left( t+1\right) ^{\nu }\right) dxdt+C\nu
^{2}\dint\limits_{Q_{T}}w_{\delta }^{2}\exp \left( 2\left( t+1\right) ^{\nu
}\right) dxdt  \label{4.21}
\end{equation}%
\begin{equation*}
-C_{1}\dint\limits_{Q_{T}}\left( \left( \nabla w_{\delta }\right)
^{2}+w_{\delta }^{2}\right) ^{2}\exp \left( 2\left( t+1\right) ^{\nu
}\right) dxdt-C\left\Vert \nabla w_{\delta }\left( x,0\right) \right\Vert
_{L_{2}\left( \Omega \right) }^{2}.
\end{equation*}%
Choose $\nu _{2}\geq \nu _{1}>1$ such that $C\sqrt{\nu _{2}}/2\geq C_{1}.$
Then, using (\ref{4.180}) and (\ref{4.21}), we obtain 
\begin{equation*}
\dint\limits_{Q_{T}}\left( w_{\delta t}-Lw_{\delta }\right) ^{2}\exp \left(
2\left( t+1\right) ^{\nu }\right) dxdt-C_{1}\dint\limits_{Q_{T}}\left(
\left( \nabla w_{\delta }\right) ^{2}+w_{\delta }^{2}\right) ^{2}\exp \left(
2\left( t+1\right) ^{\nu }\right) dxdt
\end{equation*}%
\begin{equation*}
\geq C_{1}\sqrt{\nu }\dint\limits_{Q_{T}}\left( \nabla w_{\delta }\right)
^{2}\exp \left( 2\left( t+1\right) ^{\nu }\right) dxdt+C_{1}\nu
^{2}\dint\limits_{Q_{T}}w_{\delta }^{2}\exp \left( 2\left( t+1\right) ^{\nu
}\right) dxdt
\end{equation*}%
\begin{equation}
-C_{1}\left( 1+\left\Vert v^{\ast }\right\Vert _{H^{2}\left( Q_{T}\right)
}^{2}\right)  \label{4.210}
\end{equation}%
\begin{equation*}
\geq C_{1}\sqrt{\nu }\dint\limits_{Q_{T_{\tau }}}\left( \nabla w_{\delta
}\right) ^{2}\exp \left( 2\left( t+1\right) ^{\nu }\right) dxdt+C_{1}\nu
^{2}\dint\limits_{Q_{T\tau }}w_{\delta }^{2}\exp \left( 2\left( t+1\right)
^{\nu }\right) dxdt
\end{equation*}%
\begin{equation*}
-C_{1}\left( 1+\left\Vert v^{\ast }\right\Vert _{H^{2}\left( Q_{T}\right)
}^{2}\right)
\end{equation*}%
\begin{equation*}
\geq C_{1}\exp \left( 2\left( \tau +1\right) ^{\nu }\right) \left\Vert
w_{\delta }\right\Vert _{H^{1,0}\left( Q_{T\tau }\right) }^{2}-C_{1}\left(
1+\left\Vert v^{\ast }\right\Vert _{H^{2}\left( Q_{T}\right) }^{2}\right) ,%
\text{ }\forall \nu \geq \nu _{2}.
\end{equation*}%
Hence, (\ref{4.19})-(\ref{4.210}) imply that 
\begin{equation*}
\exp \left( 2\left( T+1\right) ^{\nu }\right) \delta ^{2}\left( 1+\left\Vert
v^{\ast }\right\Vert _{H^{2}\left( Q_{T}\right) }^{2}\right) +C_{1}\left(
1+\left\Vert v^{\ast }\right\Vert _{H^{2}\left( Q_{T}\right) }^{2}\right)
\end{equation*}%
\begin{equation*}
\geq C_{1}\exp \left( 2\left( \tau +1\right) ^{\nu }\right) \left\Vert
w\right\Vert _{H^{1,0}\left( Q_{T\tau }\right) }^{2},\forall \nu \geq \nu
_{2}.
\end{equation*}%
Dividing this by $\exp \left( 2\left( \tau +1\right) ^{\nu }\right) ,$ we
obtain for all $\forall \nu \geq \nu _{2}$ 
\begin{equation}
\left\Vert w\right\Vert _{H^{1,0}\left( Q_{T\tau }\right) }^{2}\leq C_{1} 
\left[ \exp \left( 2\left( T+1\right) ^{\nu }\right) \delta ^{2}+\exp \left(
-2\left( \tau +1\right) ^{\nu }\right) \right] \left( 1+\left\Vert v^{\ast
}\right\Vert _{H^{2}\left( Q_{T}\right) }^{2}\right) .  \label{4.22}
\end{equation}%
Choose $\nu =\nu \left( \delta \right) $ such that (\ref{2}) would be
satisfied with $k=2$, i.e. 
\begin{equation}
\exp \left( 2\left( T+1\right) ^{\nu \left( \delta \right) }\right) =\frac{1%
}{\delta }.  \label{4.23}
\end{equation}%
Hence, 
\begin{equation}
\nu =\nu \left( \delta \right) =\ln \left[ \left( \ln \left( \delta
^{-1/2}\right) ^{1/\ln \left( T+1\right) }\right) \right] .  \label{4.24}
\end{equation}%
The choice (\ref{4.23}), (\ref{4.24}) is possible since (\ref{4.140}) holds
and $\delta \in \left( 0,\delta _{0}\right) .$ It follows from (\ref{2}), (%
\ref{9}) and (\ref{4.23}) that 
\begin{equation}
\exp \left( -2\left( \tau +1\right) ^{\nu \left( \delta \right) }\right)
=\exp \left[ -\frac{2}{2^{c}}\ln \left( \delta ^{-1}\right) ^{c}\right] .
\label{4.25}
\end{equation}%
Hence, using (\ref{4.22})-(\ref{4.25}) and Lemma 3.1, we obtain%
\begin{equation*}
\left\Vert w\right\Vert _{H^{1,0}\left( Q_{T\tau }\right) }^{2}\leq
C_{1}\left( 1+\left\Vert v^{\ast }\right\Vert _{H^{2}\left( Q_{T}\right)
}^{2}\right) \exp \left[ -\frac{2}{2^{c}}\ln \left( \delta ^{-1}\right) ^{c}%
\right] ,
\end{equation*}%
which implies (\ref{4.15}). \ $\square $

\section{The Global Strict Convexity}

\label{sec:5}

\subsection{The weighted Tikhonov-like functional}

\label{sec:5.1}

While the linear case (\ref{4.1}) was studied in section 4, in this section
we consider the quasilinear case. First, just like in section 4, we consider
the function $v\left( x,t\right) =u\left( x,t\right) -g\left( x\right) .$
Then, assuming (\ref{4.6}), we obtain instead of (\ref{1.6})-(\ref{1.8}):%
\begin{equation}
v_{t}=Lv+G\left( \nabla v,v,x,t\right) ,\left( x,t\right) \in Q_{T},
\label{5.1}
\end{equation}%
\begin{equation}
v\mid _{S_{T}}=0,  \label{5.2}
\end{equation}%
\begin{equation}
v\left( x,T\right) =0,x\in \Omega ,  \label{5.3}
\end{equation}%
\begin{equation}
G\left( \nabla v,v,x,t\right) =Lg+F\left( \nabla v+\nabla g,v+g,x,t\right) .
\label{5.4}
\end{equation}%
In our derivations below $v_{t},Lv$ and arguments of the function $F\left(
\nabla v+\nabla g,v+g,x,t\right) $ must be uniformly bounded for all $\left(
x,t\right) \in \overline{Q}_{T}.$ Hence, similarly with \cite{Bak,KQ}, we
now need to impose a higher smoothness than just $v\in H_{0}^{2}\left(
Q_{T}\right) $ as in section 4. Consider an integer $k$ such that $k>\left[
\left( n+1\right) /2\right] +2,$ where $\left[ \left( n+1\right) /2\right] $
denotes the maximal integer which does not exceed $\left( n+1\right) /2.$
Then embedding theorem implies that $H^{k}\left( Q_{T}\right) \subset
C^{2}\left( \overline{Q}_{T}\right) $ and 
\begin{equation}
\left\Vert f\right\Vert _{C^{2}\left( \overline{Q}_{T}\right) }\leq
E\left\Vert f\right\Vert _{H^{k}\left( Q_{T}\right) },\forall f\in
H^{k}\left( Q_{T}\right) ,  \label{5.5}
\end{equation}%
where the constant $E=E\left( Q_{T}\right) >0$ depends only on the domain $%
Q_{T}.$ Define the subspace $H_{0}^{k}\left( Q_{T}\right) \subset
H^{k}\left( Q_{T}\right) $ as%
\begin{equation*}
H_{0}^{k}\left( Q_{T}\right) =\left\{ v\in Q_{T}:v\mid _{S_{T}}=0,v\left(
x,T\right) =0\right\} .
\end{equation*}%
In addition, since we need below the function $\left( Lg\right) \left(
x,t\right) $ to be bounded in $\overline{Q}_{T},$ then we assume that 
\begin{equation}
g\in C^{2}\left( \overline{\Omega }\right) .  \label{5.50}
\end{equation}

Let $R>0$ be an arbitrary number. We consider the ball $B\left( R\right) $
in the space $H_{0}^{k}\left( Q_{T}\right) ,$ 
\begin{equation}
B\left( R\right) =\left\{ v\in H_{0}^{k}\left( Q_{T}\right) :\left\Vert
v\right\Vert _{H^{k}\left( Q_{T}\right) }<R\right\} .  \label{5.51}
\end{equation}%
Hence, by (\ref{5.5})%
\begin{equation}
\overline{B\left( R\right) }\subset C^{2}\left( \overline{Q}_{T}\right) ,
\label{5.52}
\end{equation}%
\begin{equation}
\left\Vert v\right\Vert _{C^{2}\left( \overline{Q}_{T}\right) }\leq E_{R},%
\text{ }\forall v\in \overline{B\left( R\right) },  \label{5.53}
\end{equation}%
where the number $E_{R}=E_{R}\left( Q_{T},R\right) =const.>0$ depends only
on listed parameters.

We want to find an approximate solution of problem (\ref{5.1})-(\ref{5.4})
in the closed ball $\overline{B\left( R\right) }$. To do this, we select a
number $\tau \in \left( 0,T\right) $ and minimize the following weighted
Tikhonov-like functional%
\begin{equation}
I_{\alpha ,\nu }\left( v\right) =\exp \left( -2\left( \tau +1\right) ^{\nu
}\right) \dint\limits_{Q_{T}}\left( v_{t}-Lv-G\left( \nabla v,v,x,t\right)
\right) ^{2}\exp \left( 2\left( t+1\right) ^{\nu }\right) dxdt  \label{5.6}
\end{equation}%
\begin{equation*}
+\alpha \left\Vert v\right\Vert _{H^{k}\left( Q_{T}\right) }^{2},v\in 
\overline{B\left( R\right) }.
\end{equation*}%
The multiplier $\exp \left( -2\left( \tau +1\right) ^{\nu }\right) $ in (\ref%
{5.6}) introduced to balance two terms in the right hand side of (\ref{5.6}%
). Indeed, the regularization parameter $\alpha \in \left( 0,1\right) $ and 
\begin{equation}
\min_{t\in \left[ \tau ,T\right] }\left[ \exp \left( 2\left( t+1\right)
^{\nu }\right) \exp \left( -2\left( \tau +1\right) ^{\nu }\right) \right] =1,
\label{5.06}
\end{equation}
also, see (\ref{5.7}).

\textbf{Minimization Problem 2}. \emph{Minimize the functional }$I_{\alpha
,\lambda ,\nu }\left( v\right) $\emph{\ on the ball }$\overline{B\left(
R\right) }.$

\subsection{Theorems about the functional $I_{\protect\alpha ,\protect\nu %
}\left( v\right) $}

\label{sec:5.2}

The central theorem of this section is Theorem 5.1.

\textbf{Theorem 5.1} (global strict convexity). \emph{Assume that conditions
(\ref{1.1})-(\ref{1.10}) and (\ref{5.4}) hold. Then the functional }$%
I_{\alpha ,\nu }\left( v\right) $\emph{\ has the Fr\'{e}chet derivative }$%
I_{\alpha ,\nu }^{\prime }\left( v\right) \in H_{0}^{2}\left( Q_{T}\right) $%
\emph{\ at every point }$v\in H_{0}^{2}\left( Q_{T}\right) $\emph{\ and for
all values of parameters }$\alpha ,\nu \geq 0.$\emph{\ Let }$\nu _{0}>1$%
\emph{\ be the number of Theorem 2.1.} \emph{Then there exists a number }%
\begin{equation}
\nu _{3}=\nu _{3}\left( \mu _{1},\mu _{2},\max_{i,j}\left\Vert
a_{ij}\right\Vert _{C^{1}\left( \overline{Q}_{T}\right) },\left\Vert
g\right\Vert _{C^{2}\left( \overline{\Omega }\right) },Q_{T},R,\overline{C}%
\right) \geq \nu _{0}  \label{5.60}
\end{equation}%
\emph{depending only on listed parameters such that }$2C_{2}\exp \left(
-2\left( \tau +1\right) ^{\nu }\right) \in \left( 0,1\right) $ \emph{for }$%
\nu \geq \nu _{3}$\emph{\ and if for these values of }$\nu $\emph{\ the
regularization parameter }$\alpha $\emph{\ is such that }%
\begin{equation}
\alpha \in \left[ 2C_{2}\exp \left( -2\left( \tau +1\right) ^{\nu }\right)
,1\right) ,  \label{5.61}
\end{equation}%
\emph{\ then the functional }$I_{\alpha ,\nu }\left( v\right) $\emph{\ is
strictly convex on }$\overline{B\left( R\right) }.$\emph{\ More precisely,
for every }$\tau \in \left( 0,T\right) $\emph{\ the following strict
convexity estimate holds}%
\begin{equation*}
I_{\alpha ,\nu }\left( v_{2}\right) -I_{\alpha ,\nu }\left( v_{1}\right)
-I_{\alpha ,\nu }^{\prime }\left( v_{1}\right) \left( v_{2}-v_{1}\right)
\end{equation*}%
\begin{equation}
\geq C_{2}\left\Vert h\right\Vert _{H^{1,0}\left( Q_{T\tau }\right) }^{2}+%
\frac{\alpha }{2}\left\Vert h\right\Vert _{H^{k}\left( Q_{T}\right)
}^{2},\forall v_{1},v_{2}\in \overline{B\left( R\right) },  \label{5.7}
\end{equation}%
\emph{where the constant }$C_{2}>0$\emph{\ depends only on parameters listed
in (\ref{5.60}). }

Everywhere below $C_{2}>0$ denotes different positive constants depending on
the same parameters as those listed in (\ref{5.60}).

\textbf{Remark 5.1.} \emph{The presence of the term }$C_{2}\left\Vert
h\right\Vert _{H^{1,0}\left( Q_{T\tau }\right) }^{2}$\emph{\ in the right
hand side of (\ref{5.7}) indicates that the convergence of a gradient-like
method of the minimization of functional }$I_{\alpha ,\nu }\left( v\right) $%
\emph{\ is likely faster in the space }$H^{1,0}\left( Q_{T\tau }\right) $%
\emph{\ then in the space }$H^{k}\left( Q_{T}\right) .$\emph{\ }

\textbf{Theorem 5.2}. \emph{The Fr\'{e}chet derivative }$I_{\alpha ,\nu
}^{\prime }\left( v\right) $\emph{\ of the functional }$I_{\alpha ,\nu
}\left( v\right) $\emph{\ is Lipschitz continuous on }$B\left( 2R\right) $%
\emph{\ for all values of parameters }$\alpha ,\lambda ,\nu \geq 0.$\emph{\
In other words, there exists a number }%
\begin{equation*}
D=D\left( \mu _{1},\mu _{2},\max_{i,j}\left\Vert a_{ij}\right\Vert
_{C^{1}\left( \overline{Q}_{T}\right) },\left\Vert g\right\Vert
_{H^{2}\left( \Omega \right) },Q_{T},R,\overline{C},\lambda ,\nu ,\alpha
\right) >0
\end{equation*}%
\emph{depending only on listed parameters such that }%
\begin{equation}
\left\Vert I_{\alpha ,\nu }^{\prime }\left( v_{2}\right) -I_{\alpha ,\nu
}^{\prime }\left( v_{1}\right) \right\Vert _{H^{2}\left( Q_{T}\right) }\leq
D\left\Vert v_{2}-v_{1}\right\Vert _{H^{2}\left( Q_{T}\right) },\text{ }%
\forall v_{1},v_{2}\in B\left( 2R\right) .  \label{5.8}
\end{equation}%
\emph{Furthermore, let }$\nu _{3}$\emph{\ be the number of Theorem 5.1. Then}
\emph{for every pair }$\alpha >0,$ $\nu \geq \nu _{3}$\emph{\ there exists
unique minimizer }$v_{\min }\in \overline{B\left( R\right) }$\emph{\ of the
functional }$I_{\alpha ,\nu }\left( v\right) $\emph{\ on the closed ball }$%
\overline{B\left( R\right) }$\emph{\ and the following inequality holds}%
\begin{equation}
I_{\alpha ,\nu }^{\prime }\left( v_{\min }\right) \left( v_{\min }-w\right)
\leq 0,\text{ }\forall w\in \overline{B\left( R\right) }.  \label{5.80}
\end{equation}

Let $P_{\overline{B}}:H_{0}^{k}\left( Q_{T}\right) \rightarrow \overline{%
B\left( R\right) }$ be the orthogonal projection operator mapping the space $%
H_{0}^{k}\left( Q_{T}\right) $ onto the closed ball $\overline{B\left(
R\right) }.$ Let $v_{0}\in B\left( R\right) $ be an arbitrary point of $%
B\left( R\right) .$ Let the number $\gamma \in \left( 0,1\right) .$ Consider
the sequence of the gradient projection method, 
\begin{equation}
v_{n}=P_{\overline{B}}\left( v_{n-1}-\gamma I_{\alpha ,\nu }^{\prime }\left(
v_{n-1}\right) \right) ,n=1,2,...  \label{5.9}
\end{equation}

\textbf{Theorem 5.3.} \emph{Let }$\nu _{3}$\emph{\ be the number of Theorem
5.1. Choose the number }$\nu \geq \nu _{3}.$\emph{\ Let }$v_{\min }\in 
\overline{B\left( R\right) }$\emph{\ be the unique minimizer of the
functional }$I_{\alpha ,\nu }\left( v\right) $\emph{\ on the set }$\overline{%
B\left( R\right) }$\emph{\ (Theorem 5.2). Then there exists a sufficiently
small number }%
\begin{equation*}
\gamma _{0}=\gamma _{0}\left( \mu _{1},\mu _{2},\max_{i,j}\left\Vert
a_{ij}\right\Vert _{C^{1}\left( \overline{Q}_{T}\right) },\left\Vert
g\right\Vert _{C^{2}\left( \overline{\Omega }\right) },Q_{T},R,\overline{C}%
,\nu ,\alpha \right) \in \left( 0,1\right)
\end{equation*}%
\emph{such that for every }$\gamma \in \left( 0,\gamma _{0}\right) $\emph{\
there exists a number }$\theta =\theta \left( \gamma \right) $\emph{\ such
that }%
\begin{equation}
\left\Vert v_{n}-v_{\min }\right\Vert _{H^{k}\left( Q_{T}\right) }\leq
\theta ^{n}\left\Vert v_{\min }-v_{0}\right\Vert _{H^{k}\left( Q_{T}\right)
}.  \label{5.10}
\end{equation}

Consider now the case of noise in the data. Following one of the Tikhonov's
concept of the regularization \cite{BK,T}, we assume that there exists exact
noiseless data $g^{\ast }\in H^{2}\left( \Omega \right) $ in (\ref{5.4})
and, respectively, there exists exact solution $v^{\ast }\in B\left(
R\right) $ of problem (\ref{5.1})-(\ref{5.4}). Let $\delta \in \left(
0,1\right) $ be the level of noise in the data $g\left( x\right) $, i.e.%
\begin{equation}
\left\Vert g-g^{\ast }\right\Vert _{C^{2}\left( \overline{\Omega }\right)
}<\delta .  \label{5.11}
\end{equation}%
In Theorem 5.4 we estimate the accuracy of the minimizer, i.e. the norm $%
\left\Vert v_{\min }-v^{\ast }\right\Vert _{H^{1,0}\left( Q_{T\tau }\right)
} $ for any $\tau \in \left( 0,T\right) .$ In turn, this estimate, combined
with (\ref{5.10}), provides the convergence rate of the sequence (\ref{5.9})
to the exact solution. Note that since $\delta \in \left( 0,1\right) ,$ then
by (\ref{5.11}), we replace below dependencies on $\left\Vert g\right\Vert
_{C^{2}\left( \overline{\Omega }\right) }$ of the above numbers $\nu
_{3},C_{2},D,\gamma _{0}$ with dependencies on $\left\Vert g^{\ast
}\right\Vert _{C^{2}\left( \overline{\Omega }\right) }.$

\textbf{Theorem 5.4 }(estimates of the accuracy and the convergence rate). 
\emph{Assume that the exact solution of problem (\ref{5.1})-(\ref{5.5}) }$%
v^{\ast }\in B\left( R\right) $\emph{\ and that (\ref{5.11}) holds. Let }$%
\nu _{3}>1$ \emph{be the number of Theorem 5.1. Select an arbitrary number }$%
\tau \in \left( 0,T\right) .$ \emph{For any }$\delta \in \left( 0,1\right) $%
\emph{\ set the number }$\nu =\nu \left( \delta \right) $\emph{\ be the same
as in (\ref{4.24}). Let the number }$\delta _{0}=\delta _{0}\left( \mu
_{1},\mu _{2},\max_{i,j}\left\Vert a_{ij}\right\Vert _{C^{1}\left( \overline{%
Q}_{T}\right) },\left\Vert g^{\ast }\right\Vert _{C^{2}\left( \overline{%
\Omega }\right) },Q_{T},R,\overline{C}\right) >0$\emph{\ be so small that }%
\begin{equation}
\nu \left( \delta _{0}\right) \geq \nu _{3}\text{ \emph{and} }2C_{2}\exp
\left( -2\left( \tau +1\right) ^{\nu \left( \delta _{0}\right) }\right) \in
\left( 0,1\right) .  \label{6.2}
\end{equation}%
\emph{Let} $\delta \in \left( 0,\delta _{0}\right) $ \emph{and let }$v_{\min
}\in \overline{B\left( R\right) }$\emph{\ be the unique minimizer of the
functional }$I_{\alpha ,\nu }\left( v\right) $\emph{\ on the set }$\overline{%
B\left( R\right) }$\emph{\ (Theorem 5.2). Let }$\gamma _{0}$ \emph{be the
number defined in Theorem 5.3. Let }$\gamma \in \left( 0,\gamma _{0}\right) $%
\emph{\ and }$\theta =\theta \left( \gamma \right) \in \left( 0,1\right) $%
\emph{\ also be the numbers of Theorem 5.3.\ Choose the regularization
parameter }$\alpha $ \emph{as }%
\begin{equation}
\alpha =\alpha \left( \delta \right) =2C_{2}\exp \left( -2\left( \tau
+1\right) ^{\nu \left( \delta \right) }\right) .  \label{6.3}
\end{equation}%
\emph{\ } \emph{Then the following accuracy and convergence estimates hold}%
\begin{equation}
\left\Vert v^{\ast }-v_{\min }\right\Vert _{H^{1,0}\left( Q_{T\tau }\right)
}\leq C_{2}\exp \left[ -\frac{1}{2^{c}}\ln \left( \delta ^{-1}\right) ^{c}%
\right] ,  \label{5.110}
\end{equation}%
\begin{equation}
\left\Vert v^{\ast }-v_{n}\right\Vert _{H^{1,0}\left( Q_{T\tau }\right)
}\leq C_{2}\exp \left[ -\frac{1}{2^{c}}\ln \left( \delta ^{-1}\right) ^{c}%
\right] +\theta ^{n}\left\Vert v_{\min }-v_{0}\right\Vert _{H^{k}\left(
Q_{T}\right) },  \label{5.111}
\end{equation}%
\emph{where the number }$c=c\left( \tau ,T\right) \in \left( 0,1\right) $%
\emph{\ is defined in (\ref{8}).}

\textbf{Remark 5.2}. \emph{According to section 1, since }$R>0$\emph{\ is an
arbitrary number and since the starting point }$v_{0}\in B\left( R\right) $%
\emph{\ of the gradient projection method is an arbitrary point of }$B\left(
R\right) ,$\emph{\ then Theorem 5.4 implies the \textbf{global convergence}
to the exact solution of the gradient projection method of the minimization
of the functional }$I_{\alpha ,\nu }\left( v\right) .$\emph{\ }

In this paragraph, we temporary assume that Theorem 5.1 is proved. Then the
proof of the Lipschitz continuity (\ref{5.8}) of the Fr\'{e}chet derivative $%
I_{\alpha ,\nu }^{\prime }$ is very similar to the proof of Theorem 3.1 of 
\cite{Bak}. The rest of Theorem 5.2 follows from (\ref{5.8}) and Lemma 2.1
of \cite{Bak}. Given Theorem 5.1, Theorem 5.3 follows from Theorem 3.3 of 
\cite{Bak}.

Therefore, we prove in section 6 only Theorems 5.1 and 5.4. Below, if we say
that a vector function belongs to a certain Banah space, then this means
that each of its components belongs to that space. The norm of that vector
function in that space is defined as the square root of the sum of squared
norms of its components.

\section{Proofs of Theorems 5.1 and 5.4}

\label{sec:6}

\subsection{Proof of Theorem 5.1}

\label{sec:6.1}

Consider two arbitrary functions $v_{1},v_{2}\in B\left( R\right) .$ Denote $%
h=v_{2}-v_{1}.$ Then 
\begin{equation}
h\in B\left( 2R\right) .  \label{5.12}
\end{equation}%
Hence, by (\ref{5.53}) 
\begin{equation}
h\in C^{2}\left( \overline{Q}_{T}\right) ,\left\Vert h\right\Vert
_{C^{2}\left( \overline{Q}_{T}\right) }\leq E_{2R}.  \label{5.13}
\end{equation}%
Consider the expression 
\begin{equation*}
v_{2t}-Lv_{2}-G\left( \nabla v_{2},v_{2},x,t\right)
\end{equation*}%
\begin{equation}
=\left( h_{t}-Lh\right) -G\left( \nabla v_{1}+\nabla h,v_{1}+h,x,t\right)
+\left( v_{1t}-Lv_{1}\right) .  \label{5.14}
\end{equation}%
Using the multidimensional analog of Taylor formula \cite{V}, (\ref{1.10}), (%
\ref{5.4}), (\ref{5.53}) and (\ref{5.13}), we obtain%
\begin{equation}
G\left( \nabla v_{1}+\nabla h,v_{1}+h,x,t\right)  \label{5.15}
\end{equation}%
\begin{equation*}
=G\left( \nabla v_{1},v_{1},x,t\right) +G_{1}\left( x,t\right) \nabla
h+G_{2}\left( x,t\right) h+G_{3}\left( x,t,\nabla h,h\right) ,
\end{equation*}%
where the $n-$D vector function $G_{1}\left( x,t\right) \in C\left( 
\overline{Q}_{T}\right) ,$ the function $G_{2}\left( x,t\right) \in C\left( 
\overline{Q}_{T}\right) $ and the function $G_{3}\left( x,t,\nabla
h,h\right) $ is such that%
\begin{equation}
\left\vert G_{3}\left( x,t,\nabla h,h\right) \right\vert \leq C_{2}\left(
\left( \nabla h\right) ^{2}+h^{2}\right) \left( x,t\right) ,\forall h\in
B\left( 2R\right) ,\forall \left( x,t\right) \in \overline{Q}_{T}.
\label{5.16}
\end{equation}%
In addition, 
\begin{equation}
\left\Vert G_{1}\left( x,t\right) \right\Vert _{C\left( \overline{Q}%
_{T}\right) },\left\Vert G_{2}\left( x,t\right) \right\Vert _{C\left( 
\overline{Q}_{T}\right) }\leq C_{2}.  \label{5.17}
\end{equation}

Hence, (\ref{5.14}) and (\ref{5.15}) imply%
\begin{equation*}
\left( v_{2t}-Lv_{2}-G\left( \nabla v_{2},v_{2},x,t\right) \right)
^{2}-\left( v_{1t}-Lv_{1}-G\left( \nabla v_{1},v_{1},x,t\right) \right) ^{2}
\end{equation*}%
\begin{equation}
=2\left( v_{1t}-Lv_{1}\right) \left[ \left( h_{t}-Lh\right) -G_{1}\left(
x,t\right) \nabla h-G_{2}\left( x,t\right) h\right]  \label{5.18}
\end{equation}%
\begin{equation*}
+\left[ \left( h_{t}-Lh\right) -G_{1}\left( x,t\right) \nabla h-G_{2}\left(
x,t\right) h-G_{3}\left( x,t,\nabla h,h\right) \right] ^{2}
\end{equation*}%
\begin{equation*}
-2\left( v_{1t}-Lv_{1}\right) G_{3}\left( x,t,\nabla h,h\right) =Lin\left(
h\right) \left( x,t\right) +Nonlin\left( h\right) \left( x,t\right) ,
\end{equation*}%
where $Lin\left( h\right) $ and $Nonlin\left( h\right) $ denote linear and
nonlinear expressions with respect to $h$ respectively,%
\begin{equation}
Lin\left( h\right) \left( x,t\right) =2\left( v_{1t}-Lv_{1}\right) \left[
\left( h_{t}-Lh\right) -G_{1}\left( x,t\right) \nabla h-G_{2}\left(
x,t\right) h\right] .  \label{5.19}
\end{equation}%
We represent the term $Nonlin\left( h\right) \left( x,t\right) $ in (\ref%
{5.18}) as%
\begin{equation*}
Nonlin\left( h\right) \left( x,t\right) =\left( h_{t}-Lh\right) ^{2}
\end{equation*}%
\begin{equation*}
-2\left( h_{t}-Lh\right) \left[ G_{1}\left( x,t\right) \nabla h+G_{2}\left(
x,t\right) h+G_{3}\left( x,t,\nabla h,h\right) \right]
\end{equation*}%
\begin{equation*}
+\left[ G_{1}\left( x,t\right) \nabla h+G_{2}\left( x,t\right) h+G_{3}\left(
x,t,\nabla h,h\right) \right] ^{2}.
\end{equation*}%
Hence, using Cauchy-Schwarz inequality, (\ref{5.16}) and (\ref{5.17}), we
obtain%
\begin{equation}
Nonlin\left( h\right) \left( x,t\right) \geq \frac{1}{2}\left(
h_{t}-Lh\right) ^{2}-C_{2}\left( \left( \nabla h\right) ^{2}+h^{2}\right)
,\forall \left( x,t\right) \in Q_{T}.  \label{5.20}
\end{equation}%
Also, by (\ref{5.5}) and (\ref{5.16})-(\ref{5.18}) 
\begin{equation}
\left\vert Nonlin\left( h\right) \right\vert \left( x,t\right) \leq
C_{2}\left\Vert h\right\Vert _{C^{2}\left( \overline{Q}_{T}\right) }^{2}\leq
C_{2}\left\Vert h\right\Vert _{H^{k}\left( Q_{T}\right) }^{2},\forall \left(
x,t\right) \in Q_{T}.  \label{5.21}
\end{equation}

Thus, (\ref{5.6}), (\ref{5.18}) and (\ref{5.19}) imply that%
\begin{equation*}
I_{\alpha ,\nu }\left( v_{1}+h\right) -I_{\alpha ,\nu }\left( v_{1}\right)
\end{equation*}%
\begin{equation}
=\exp \left( -2\left( \tau +1\right) ^{\nu }\right)
\dint\limits_{Q_{T}}Lin\left( h\right) \exp \left( 2\left( t+1\right) ^{\nu
}\right) dxdt+2\alpha \left( v_{1},h\right) _{k}  \label{5.22}
\end{equation}%
\begin{equation*}
+\exp \left( -2\left( \tau +1\right) ^{\nu }\right)
\dint\limits_{Q_{T}}Nonlin\left( h\right) \exp \left( 2\left( t+1\right)
^{\nu }\right) dxdt+\alpha \left\Vert h\right\Vert _{H^{k}\left(
Q_{T}\right) }^{2},
\end{equation*}%
where $\left( .,.\right) _{k}$ is the scalar product in $H^{k}\left(
Q_{T}\right) .$ Assuming temporary that $h$ is an arbitrary function of $%
H^{k}\left( Q_{T}\right) ,$ consider the expression $X\left( h\right) ,$%
\begin{equation}
X\left( h\right) =\exp \left( -2\left( \tau +1\right) ^{\nu }\right)
\dint\limits_{Q_{T}}Lin\left( h\right) \exp \left( 2\left( t+1\right) ^{\nu
}\right) dxdt+2\alpha \left( v_{1},h\right) _{k}.  \label{5.23}
\end{equation}%
It follows from (\ref{5.53}), (\ref{5.17}) and (\ref{5.19}) that $%
X:H^{k}\left( Q_{T}\right) \rightarrow \mathbb{R}$ is a bounded linear
functional. Hence, by Riesz theorem there exists a unique function $%
\widetilde{X}\in H^{k}\left( Q_{T}\right) $ such that 
\begin{equation}
X\left( h\right) =\left( \widetilde{X},h\right) _{k}.  \label{5.24}
\end{equation}%
At the same time, it follows from (\ref{5.21})-(\ref{5.24}) that%
\begin{equation}
I_{\alpha ,\nu }\left( v_{1}+h\right) -I_{\alpha ,\nu }\left( v_{1}\right)
=\left( \widetilde{X},h\right) _{k}+o\left( \left\Vert h\right\Vert
_{H^{k}\left( Q_{T}\right) }\right) ,  \label{5.25}
\end{equation}%
\begin{equation*}
\lim_{\left\Vert h\right\Vert _{H^{k}\left( Q_{T}\right) }\rightarrow
0}\left( \frac{o\left( \left\Vert h\right\Vert _{H^{k}\left( Q_{T}\right)
}\right) }{\left\Vert h\right\Vert _{H^{k}\left( Q_{T}\right) }}\right) =0.
\end{equation*}%
Hence, $\widetilde{X}$ is the Fr\'{e}chet derivative $I_{\alpha ,\nu
}^{\prime }\left( v_{1}\right) \in H^{k}\left( Q_{T}\right) $ of the
functional $I_{\alpha ,\nu }\left( v\right) $ at the point $v_{1},$%
\begin{equation}
\left( \widetilde{X},h\right) _{k}=I_{\alpha ,\nu }^{\prime }\left(
v_{1}\right) \left( h\right) ,\text{ }\forall h\in H^{k}\left( Q_{T}\right) .
\label{5.26}
\end{equation}

We now come back again to the case when $h\in B\left( 2R\right) $ as in (\ref%
{5.12}). Using (\ref{5.20}), (\ref{5.22}) and (\ref{5.24})-(\ref{5.26}), we
obtain%
\begin{equation*}
I_{\alpha ,\nu }\left( v_{1}+h\right) -I_{\alpha ,\nu }\left( v_{1}\right)
-I_{\alpha ,\nu }^{\prime }\left( v_{1}\right) \left( h\right)
\end{equation*}%
\begin{equation}
\geq \frac{1}{2}\exp \left( -2\left( \tau +1\right) ^{\nu }\right)
\dint\limits_{Q_{T}}\left( h_{t}-Lh\right) ^{2}\exp \left( 2\left(
t+1\right) ^{\nu }\right) dxdt  \label{5.27}
\end{equation}%
\begin{equation*}
-C_{2}\exp \left( -2\left( \tau +1\right) ^{\nu }\right)
\dint\limits_{Q_{T}}\left( \left( \nabla h\right) ^{2}+h^{2}\right) \exp
\left( 2\left( t+1\right) ^{\nu }\right) dxdt+\alpha \left\Vert h\right\Vert
_{H^{k}\left( Q_{T}\right) }^{2}.
\end{equation*}%
We now use the Carleman estimate (\ref{1.9}). Recalling that $h\left(
x,T\right) =0$ and using (\ref{5.27}), we obtain%
\begin{equation*}
I_{\alpha ,\nu }\left( v_{1}+h\right) -I_{\alpha ,\nu }\left( v_{1}\right)
-I_{\alpha ,\nu }^{\prime }\left( v_{1}\right) \left( h\right)
\end{equation*}%
\begin{equation*}
\geq C\sqrt{\nu }\exp \left( -2\left( \tau +1\right) ^{\nu }\right)
\dint\limits_{Q_{T}}\left( \nabla h\right) ^{2}\exp \left( 2\left(
t+1\right) ^{\nu }\right) dxdt
\end{equation*}%
\begin{equation}
+C\nu ^{2}\exp \left( -2\left( \tau +1\right) ^{\nu }\right)
\dint\limits_{Q_{T}}h^{2}\exp \left( 2\left( t+1\right) ^{\nu }\right) dxdt
\label{5.28}
\end{equation}%
\begin{equation*}
-C_{2}\exp \left( -2\left( \tau +1\right) ^{\nu }\right)
\dint\limits_{Q_{T}}\left( \left( \nabla h\right) ^{2}+h^{2}\right) \exp
\left( 2\left( t+1\right) ^{\nu }\right) dxdt
\end{equation*}%
\begin{equation*}
-C\exp \left( -2\left( \tau +1\right) ^{\nu }\right) \left\Vert \nabla
h\left( x,0\right) \right\Vert _{L_{2}\left( \Omega \right) }^{2}+\alpha
\left\Vert h\right\Vert _{H^{k}\left( Q_{T}\right) }^{2}.
\end{equation*}%
By (\ref{5.5}) 
\begin{equation}
\left\Vert \nabla h\left( x,0\right) \right\Vert _{L_{2}\left( \Omega
\right) }^{2}\leq C_{2}\left\Vert h\right\Vert _{H^{k}\left( Q_{T}\right)
}^{2}.  \label{5.280}
\end{equation}

Choose the number $\nu _{3}\geq \nu _{0}>1$ depending on the same parameters
as those listed in (\ref{5.60}) and such that $C\sqrt{\nu _{3}}/2>C_{2}.$
Hence, (\ref{5.06}), (\ref{5.27})-(\ref{5.280}) and (\ref{5.61}) imply that
for all $\nu \geq \nu _{3}$ 
\begin{equation*}
I_{\alpha ,\nu }\left( v_{1}+h\right) -I_{\alpha ,\nu }\left( v_{1}\right)
-I_{\alpha ,\nu }^{\prime }\left( v_{1}\right) \left( h\right)
\end{equation*}%
\begin{equation*}
\geq C_{2}\sqrt{\nu }\exp \left( -2\left( \tau +1\right) ^{\nu }\right)
\dint\limits_{Q_{T}}\left( \nabla h\right) ^{2}\exp \left( 2\left(
t+1\right) ^{\nu }\right) dxdt
\end{equation*}%
\begin{equation*}
+C_{2}\nu ^{2}\exp \left( -2\left( \tau +1\right) ^{\nu }\right)
\dint\limits_{Q_{T}}h^{2}\exp \left( 2\left( t+1\right) ^{\nu }\right) dxdt
\end{equation*}%
\begin{equation*}
+\left\Vert h\right\Vert _{H^{k}\left( Q_{T}\right) }^{2}\left( \alpha
-C_{2}\exp \left( -2\left( \tau +1\right) ^{\nu }\right) \right)
\end{equation*}%
\begin{equation*}
\geq C_{2}\sqrt{\nu }\exp \left( -2\left( \tau +1\right) ^{\nu }\right)
\dint\limits_{Q_{T\tau }}\left( \nabla h\right) ^{2}\exp \left( 2\left(
t+1\right) ^{\nu }\right) dxdt
\end{equation*}%
\begin{equation*}
+C_{2}\nu ^{2}\exp \left( -2\left( \tau +1\right) ^{\nu }\right)
\dint\limits_{Q_{T\tau }}h^{2}\exp \left( 2\left( t+1\right) ^{\nu }\right)
dxdt+\frac{\alpha }{2}\left\Vert h\right\Vert _{H^{k}\left( Q_{T}\right)
}^{2}
\end{equation*}%
\begin{equation*}
\geq C_{2}\left\Vert h\right\Vert _{H^{1,0}\left( Q_{T\tau }\right) }^{2}+%
\frac{\alpha }{2}\left\Vert h\right\Vert _{H^{k}\left( Q_{T}\right) }^{2}.%
\text{ \ }\square
\end{equation*}

\subsection{Proof of Theorem 5.4}

\label{sec:6.2}

Denote%
\begin{equation}
I_{\alpha ,\nu }^{0}\left( v\right) =\exp \left( -2\left( \tau +1\right)
^{\nu }\right) \dint\limits_{Q_{T}}\left( v_{t}-Lv-G\left( \nabla
v,v,x,t\right) \right) ^{2}\exp \left( 2\left( t+1\right) ^{\nu }\right)
dxdt.  \label{5.29}
\end{equation}%
By (\ref{5.6}) and (\ref{5.29})%
\begin{equation}
I_{\alpha ,\nu }\left( v\right) =I_{\alpha ,\nu }^{0}\left( v\right) +\alpha
\left\Vert v\right\Vert _{H^{k}\left( Q_{T}\right) }^{2}.  \label{5.30}
\end{equation}%
By (\ref{5.4}) 
\begin{equation*}
G\left( \nabla v,v,x,t\right) =Lg+F\left( \nabla v+\nabla g,v+g,x,t\right)
\end{equation*}%
\begin{equation*}
=Lg^{\ast }+\left( Lg-Lg^{\ast }\right) +F\left[ \left( \nabla v+\nabla
g^{\ast }\right) +\left( \nabla g-\nabla g^{\ast }\right) ,\left( v+g^{\ast
}\right) +\left( g-g^{\ast }\right) ,x,t\right] .
\end{equation*}%
Hence, using the multidimensional analog of Taylor formula \cite{V} and (\ref%
{5.11}), we obtain similarly with (\ref{5.15})-(\ref{5.17})%
\begin{equation}
G\left( \nabla v,v,g,x,t\right) =Lg^{\ast }+F\left( \nabla v^{\ast }+\nabla
g^{\ast },v^{\ast }+g^{\ast },x,t\right) +P\left( x,t\right)  \label{5.31}
\end{equation}%
\begin{equation*}
=G\left( \nabla v^{\ast },v^{\ast },g^{\ast },x,t\right) +P\left( x,t\right)
.
\end{equation*}%
where the function $P\left( x,t\right) $ is such that $\left\Vert
P\right\Vert _{L_{2}\left( Q_{T}\right) }\leq C_{2}\delta .$ Since $%
v_{t}^{\ast }-Lv^{\ast }-G\left( \nabla v^{\ast },v^{\ast },g^{\ast
},x,t\right) =0,$ then%
\begin{equation*}
I_{\alpha ,\nu }^{0}\left( v^{\ast }\right) =\exp \left( -2\left( \tau
+1\right) ^{\nu }\right) \times
\end{equation*}%
\begin{equation*}
\times \dint\limits_{Q_{T}}\left[ v_{t}^{\ast }-Lv^{\ast }-G\left( \nabla
v^{\ast },v^{\ast },g^{\ast },x,t\right) +P\left( x,t\right) \right]
^{2}\exp \left( 2\left( t+1\right) ^{\nu }\right) dxdt
\end{equation*}%
\begin{equation*}
=\exp \left( -2\left( \tau +1\right) ^{\nu }\right)
\dint\limits_{Q_{T}}P^{2}\left( x,t\right) \exp \left( 2\left( t+1\right)
^{\nu }\right) dxdt
\end{equation*}%
\begin{equation*}
\leq C_{2}\exp \left( 2\left( T+1\right) ^{\nu }\right) \delta ^{2}.
\end{equation*}%
Hence, using (\ref{5.30}), we obtain%
\begin{equation}
I_{\alpha ,\nu }\left( v^{\ast }\right) \leq C_{2}\left( \exp \left( 2\left(
T+1\right) ^{\nu }\right) \delta ^{2}+\alpha \right) .  \label{5.32}
\end{equation}

Next, by (\ref{5.7}) 
\begin{equation*}
I_{\alpha ,\nu }\left( v^{\ast }\right) -I_{\alpha ,\nu }\left( v_{\min
}\right) -I_{\alpha ,\nu }^{\prime }\left( v_{\min }\right) \left( v^{\ast
}-v_{\min }\right) 
\end{equation*}%
\begin{equation}
\geq C_{2}\left\Vert v^{\ast }-v_{\min }\right\Vert _{H^{1,0}\left( Q_{T\tau
}\right) }^{2}+\frac{\alpha }{2}\left\Vert v^{\ast }-v_{\min }\right\Vert
_{H^{k}\left( Q_{T}\right) }^{2}.  \label{5.34}
\end{equation}%
By (\ref{5.80}) $-I_{\alpha ,\lambda ,\nu }^{\prime }\left( v_{\min }\right)
\left( v^{\ast }-v_{\min }\right) \leq 0.$ Hence, (\ref{5.32}) and (\ref%
{5.34}) imply that 
\begin{equation}
\left\Vert v^{\ast }-v_{\min }\right\Vert _{H^{1,0}\left( Q_{T\tau }\right)
}^{2}\leq C_{2}\left( \exp \left( 2\left( T+1\right) ^{\nu }\right) \delta
^{2}+\alpha \right) .  \label{5.35}
\end{equation}%
Recall that the numbers $\nu =\nu \left( \delta \right) $ and $\alpha
=\alpha \left( \delta \right) $ are defined in (\ref{4.24}) and (\ref{6.3})
respectively. These choices of $\nu \left( \delta \right) $ and $\alpha
\left( \delta \right) $ are possible since (\ref{6.2}) holds and $\delta \in
\left( 0,\delta _{0}\right) .$ Thus, condition (\ref{5.61}) of Theorem 5.1
imposed on $\alpha $ is in place$.$ Hence, (\ref{4.23}) and (\ref{4.25})
hold. Hence, using (\ref{3}), (\ref{4.23})-(\ref{4.25}), and (\ref{5.35}),
we obtain%
\begin{equation*}
\left\Vert v^{\ast }-v_{\min }\right\Vert _{H^{1,0}\left( Q_{T\tau }\right)
}^{2}\leq C_{2}\left( \delta +\exp \left[ -\frac{2}{2^{c}}\ln \left( \delta
^{-1}\right) ^{c}\right] \right) \leq C_{2}\exp \left[ -\frac{2}{2^{c}}\ln
\left( \delta ^{-1}\right) ^{c}\right] ,
\end{equation*}%
which implies (\ref{5.110}). To prove (\ref{5.111}), we use the triangle
inequality,%
\begin{equation}
\left\Vert v^{\ast }-v_{n}\right\Vert _{H^{1,0}\left( Q_{T\tau }\right)
}\leq \left\Vert v^{\ast }-v_{\min }\right\Vert _{H^{1,0}\left( Q_{T\tau
}\right) }+\left\Vert v_{\min }-v_{n}\right\Vert _{H^{1,0}\left( Q_{T\tau
}\right) }  \label{5.36}
\end{equation}%
\begin{equation*}
\leq \left\Vert v^{\ast }-v_{\min }\right\Vert _{H^{1,0}\left( Q_{T\tau
}\right) }+\left\Vert v_{\min }-v_{n}\right\Vert _{H^{k}\left( Q_{T}\right)
}.
\end{equation*}%
Using (\ref{5.10}), (\ref{5.110}) and (\ref{5.36}), we obtain (\ref{5.111}). 
$\square $

\begin{center}
\textbf{Acknowledgment}
\end{center}

The work of M.V. Klibanov was supported by US Army Research Laboratory and
US Army Research Office grant W911NF-19-1-0044.

\end{document}